\documentclass[preprint,12pt]{elsarticle}
\usepackage{filecontents}

\usepackage{natbib}
\usepackage{pslatex}
\usepackage{hyperref}
\usepackage{tikz}
\usepackage{amssymb}
\usepackage{graphicx}
\usepackage{marvosym}
\usetikzlibrary{patterns}
\usetikzlibrary{plotmarks}
\usetikzlibrary{arrows,snakes,backgrounds}
\def\prova{{\noindent{\textit{Proof} }}}
\newtheorem{tdef}{Definition}
\newtheorem{tteo}{Theorem}
\newtheorem{tcor}[tteo]{Corollary}
\newtheorem{tlem}{Lemma}
\newtheorem{tcon}{Condition}
\newcommand{\R}{\mathbf{R}}

\journal{ArXiv.org}

\begin{document}

\begin{frontmatter}

\title{\Large{\bf Local Pareto optimality conditions for vector quadratic fractional optimization problems}}

\author{ {\bf W. A. Oliveira \footnote{Corresponding author}}}
\address{School of Applied Sciences, University of Campinas, R. Pedro Zaccaria, 1300, Limeira 13484-350, S\~ao Paulo, Brazil. E-mail: washington.oliveira@fca.unicamp.br} 

\author{ {\bf M. A. Rojas-Medar}}
\address{Instituto de Alta Investigaci\'on, Universidad de Tarapac\'a, Casilla 7D, Arica, Chile. E-mail: marko.medar@gmail.com}

\author{ {\bf A. Beato-Moreno}}
 \address{Department of Statistics and Operations Research, College of Mathematics, University of Sevilla, 41012, Sevilla, Spain. Tel.: +34-95-4557942. E-mail: beato@us.es} 

\author{ {\bf M. B. Hern\'andez-Jim\'enez}}
\address{Departamento de Econom\'ia, M\'etodos cuantitativos e H. Econ\'omica, Universidad Pablo de Olavide, Sevilla, Spain. E-mail: mbherjim@upo.es} 

\begin{abstract}
There are several concepts and definitions that characterize and give optimality conditions for solutions of a vector optimization problem. One of the most important is 
the first-order necessary optimality condition that generalizes the Karush-Kuhn-Tucker condition. This condition ensures the existence of an arbitrary neighborhood that 
contains an local optimal solution. The present work we introduce an alternative concept to identify the local optimal solution neighborhood in vector optimization problems. 
The main aspect of this contribution is the development of necessary and sufficient Pareto optimality conditions for the solutions of a particular vector 
optimization problem, where each objective function consists of a ratio quadratic functions and the feasible set is defined by linear inequalities. We show how to 
calculate the largest radius of the spherical region centered on a local Pareto solution in which this solution is optimal. In this process we may conclude 
that the solution is also globally optimal. These conditions might be useful to determine termination criteria in the development of algorithms, 
including more general problems in which quadratic approximations are used locally.
\end{abstract}

 \begin{keyword}
Pareto optimality conditions ;  efficient solutions ; vector optimization ;  vector quadratic fractional optimization problems 
 \end{keyword}

\end{frontmatter}

\section{Introduction}\label{sec:introducao}

In general the vector optimization problem (VOP) appears in the processes for decision-making and is represented as the following problem:

\begin{center}
\begin{tabular}{lll}
\textbf{(VOP)}&\quad \mbox{\textbf{Minimize}}\quad &$\ f(x)=\left( f_1 (x),\ldots,f_m (x) \right)$\\
&\quad \mbox{subject to}&\ $h_j(x)\leqq 0,\quad j \in J$,\\
&&\ $x\in \Omega \subseteq \R^n$,
\end{tabular}
\end{center}
where $f_i: \Omega\subseteq \R^n \longrightarrow \R^m$, $i\in I\equiv \{1,\ldots,m\}$, are functions of $n$ real variables and $\Omega$ is a nonempty 
subset. The fields of real functions $h_j's$ contains $\Omega$. We denote by $S$ the \textit{feasible set} 
that is the intersection of $\Omega$ with the set of points $x$ in which $h_j(x) \leqq 0$, $j\in J\equiv\{1,\ldots,p\}$. If $x\in S$, we say 
that $x$ is a \textit{feasible point}. $f_i(x)$ is the result of the $i^{th}$ objective function if the decision maker chooses the 
action $x\in S$.  

Let $\mathbf{R}_{+}$ denote the nonnegative real numbers and $x^T$ denote the transpose of the vector $x\in \mathbf{R}^n$. Furthermore, we will 
adopt the following conventions for inequalities among vectors. If $x=(x_1,\ldots,x_m)^T \in \mathbf{R}^m$ and $y=(y_1,\ldots,y_m)^T \in \mathbf{R}^m$, then 
\begin{itemize}
\item[] $x=y$\quad if and only if\quad $x_i = y_i$,\quad for all $i\in I$;
\item[] $x<y$\quad if and only if\quad $x_i < y_i$,\quad for all $i\in I$; 
\item[] $x\leqq y$\quad if and only if\quad $x_i \leqq y_i$,\quad for all $i\in I$;
\item[] $x\leq y$\quad if and only if\quad $x\leqq y$ and $x\neq y$. 
\end{itemize}
Similarly we consider the equivalent convention for inequalities $>, \geqq$ and $\geq$. 

We say that $x\in S$ \textit{dominates\label{def:Dominancia}} $z\in S$ in (VOP) if $f(x)\leq f(z)$, that $N(y)\subseteq S$ is a \textit{neighborhood} of $y\in S$, 
that $B(x,r)$ is an \textit{open ball} centered on $x$, the radius $r$ and the \textit{boundary} $\partial B(x,r)$, and that $\bar{B}(x,r)$ is a \textit{closed ball}, defined by Euclidean distance. 

Different optimality notions for the problem (VOP) are referred to as a Pareto optimal solution~\cite{Livro:Pareto}, two of which are defined as follows.

{\tdef{\textcolor{white}{.}}}A feasible point $x^*$ is said to be a \textit{(locally) Pareto-efficient optimal solution} of (VOP), if there does not exist 
another ($x\in N(x^*)$) $x\in S$ such that $f(x)\leq f(x^*)$.

{\tdef{\textcolor{white}{.}}}A feasible point $x^*$ is said to be a \textit{(locally) weakly Pareto-efficient optimal solution} of (VOP), if there does not exist 
another ($x\in N(x^*)$) $x\in S$ such that $f(x)< f(x^*)$.\\

Every Pareto-efficient optimal solution is a weakly Pareto-efficient optimal solution. The Pareto-efficient optimal solution set is denoted by \textit{Eff}$(VOP)$, and 
locally Pareto-efficient optimal solution by \textit{Leff}$(VOP)$. The Pareto-efficient optimal solutions are also known as globally Pareto-efficient optimal solutions. 
Let the set \textit{Leff}$(VOP)\subseteq \R^n$, we say that $f($\textit{Leff}$(VOP))\subseteq \R^m$ is the \textit{Pareto-optimal curve}.

There are many contributions, concepts, and definitions that characterize and give the Pareto-efficient optimality conditions for solutions 
of a vector optimization problem (see, for instances~\cite{Livro:Chankong,Livro:Jahn,Livro:Miettinen,Tese:Osuna,Livro:Romero}). One of the most 
important is the first-order necessary optimality condition that generalizes the Karush-Kuhn-Tucker (KKT) condition. However, to obtain the 
sufficient optimality conditions, it is necessary to impose additional assumptions (like convexity and its generalizations) in the objective 
functions and in the constraint set. Otherwise, is only possible to characterize the locally Pareto-efficient optimal solutions. In general 
the non-equivalence between locally and globally Pareto-efficient optimal solutions is a difficult problem to be faced. The following example 
shows a simple problem which illustrates this. We recall that $\nabla f(x)$ denote the gradient of the function $f:\mathbf{R}^n\to \mathbf{R}$ at point $x$.\\
 
\noindent{Exemplo:{\hspace{0.2cm}}\label{exemp:Exemplo1}}Consider the following problem:

\begin{center}
\begin{tabular}{ll}
$Minimize$\quad &$f(x) = \left( 4x_1^2 - x_2^2,\ -(x_1-2)^2 + 4(x_2+1)^2 \right)$\\
\mbox{subject to}& $x=(x_1,x_2)\in \mathbf{R}^2,$
\end{tabular}
\end{center}
which their objective functions can be expressed as

\begin{center}
\begin{tabular}{ccccl}
$f_1(x)$&$=$&0,5 $x^T A_1 x + b_1^T x + c_1$&$=$& $4x_1^2 - x_2^2$,\\
$f_2(x)$&$=$&0,5 $x^T A_2 x + b_2^T x + c_2$&$=$& $-(x_1-2)^2 + 4(x_2+1)^2$,
\end{tabular}
\end{center}
where $c_1=0$, $c_2=0$, $b_1=(0,0)^T$, $b_2=(4,8)^T$, and whose matrices $A_1$ and $A_2$ are

\begin{center}
\begin{tabular}{ccc}
$A_1=\left( \begin{array}{cr} 8 & 0 \\ 0 & -2\end{array} \right)$&\ \ and\ \ &$A_2=\left( \begin{array}{rcc} -2 && 0 \\ 0 && 8\end{array} \right).$
\end{tabular}
\end{center}

This problem is an unconstrained two objectives optimization problem, which the objectives are quadratic functions with indefinite Hessian matrices. 
If $d\neq 0$, it can not simultaneously occurs 
\begin{eqnarray}\label{eq:CondicaoExemplo1}
d^T A_1 d \leqq 0\ \mbox{and}\ d^T A_2 d \leqq 0,\ \mbox{for all}\ d\in \partial B(0,1).
\end{eqnarray}
In fact, let $d=(d_1,d_2)^T\neq (0,0)^T$. If $d^T A_1 d = 8d_1^2 - 2d_2^2 \leqq 0$ and $d^T A_2 d = -2d_1^2 + 8d_2^2 \leqq 0$, then 
\begin{eqnarray*}
0< 6d_1^2 + 6 d_2^2 = \left(8d_1^2 - 2d_2^2\right) + \left(8d_2^2 - 2d_1^2\right) \leqq 0.
\end{eqnarray*}This is a contradiction, therefore inequalities~(\ref{eq:CondicaoExemplo1}) only occur if $d\equiv0$. A necessary condition 
for a point $x^*$ be a locally Pareto-efficient optimal solution~\cite{Livro:Chankong} is that there are real numbers $\tau_1, \tau_2 \geqq 0$, 
not all zero, such that 
\begin{eqnarray}\label{eq:CondicaoNecessariaExemplo1}
\tau_1\nabla f_1(x^*) + \tau_2\nabla f_2(x^*)=0.
\end{eqnarray}

If $\tau_2=0$ and $\tau_1>0$, then $x^*=(0,0)^T$. In the direction $d=(\frac{1}{\sqrt{5}},\frac{-2}{\sqrt{5}})^T\in \partial B(0,1)$, for 
every point of the form $x'=x^* +\lambda d$, for each real value $\lambda\in(0,\frac{4\sqrt{5}}{5})$, 
simultaneously occur $f_1(x')=f_1(x^*)$ and $f_2(x')<f_2(x^*)$. Therefore, the point $x^*=(0,0)^T$ is not a locally Pareto-efficient optimal solution. 

If $\tau_1=0$ and $\tau_2>0$, then $x^*=(2,-1)^T$. In the direction $d=(\frac{-2}{\sqrt{5}},\frac{1}{\sqrt{5}})^T\in \partial B(0,1)$, for 
every point of the form $x'=x^* +\lambda d$, for each real value $\lambda\in(0,2\sqrt{5})$, 
simultaneously occurs $f_1(x')<f_1(x^*)$ and $f_2(x')=f_2(x^*)$. Therefore, the point $x^*=(2,-1)^T$ is not a locally Pareto-efficient optimal 
solution. 

Now, if $\tau_1, \tau_2 >0$, an equivalent way to write the condition~(\ref{eq:CondicaoNecessariaExemplo1}) is 
\begin{eqnarray}\label{eq:CondicaoSuficienteExemplo1}
\tau\nabla f_1(x^*) + \nabla f_2(x^*)=0,\ \mbox{where}\ \tau=\frac{\tau_1}{\tau_2} > 0.
\end{eqnarray}
But, this condition is also sufficient for the points $x^*$, unlike $(0,0)^T$ and $(2,-1)^T$, are locally Pareto-efficient optimal solution. In fact, 
suppose that $x^*$ is not locally Pareto-efficient optimal solution, then there exists a direction $d\in \partial B(0,1)$ and a number 
real $\lambda>0$ such that $f(x^* +\lambda d)\leq f(x^*)$ is valid. That is, there is a descent direction $d$ which locally occurs 
$\nabla f_1(x^*)^T d\leqq 0$ and $\nabla f_2(x^*)^T d\leqq 0$, with $\nabla f_1(x^*)\neq 0$ and $\nabla f_2(x^*)\neq 0$. However, for $x^*$ 
and $\tau>0$ if the equation~(\ref{eq:CondicaoSuficienteExemplo1}) is valid, by Stiemke's alternative theorem~\cite{Art:stiemke1915positive}, 
$\nabla f_1(x^*)^T d<0$ and $\nabla f_2(x^*)^T d \leqq 0$ has no solution, and $\nabla f_1(x^*)^T d \leqq 0$ and 
$\nabla f_2(x^*)^T d < 0$ has no solution.

Thus, when $f(x^* +\lambda d)\leq f(x^*)$ and Stiemke's alternative theorem are valid to must occur $\nabla f_1(x^*)^T d = 0$ and 
$\nabla f_2(x^*)^T d = 0$, and to must simultaneously occurs $d^T A_1 d \leqq 0$ and $d^T A_2 d < 0$, or 
$d^T A_1 d < 0$ and $d^T A_2 d \leqq 0$. What is impossible, because~(\ref{eq:CondicaoExemplo1}). Therefore, the locally 
Pareto-efficient optimal solutions are those that satisfy~(\ref{eq:CondicaoSuficienteExemplo1}), that is, $\displaystyle 8\tau x_1^* - 2(x_1^* - 2) =0$ and $\displaystyle -2\tau x_2^* + 8(x_2^* + 1) =0$, 
that is,
\begin{eqnarray}\label{eq:PontosEficientesExemplo1}
\left\{\begin{array}{cll} x_1^* = & \frac{2}{-4\tau + 1} &,\ 0< \tau \neq \frac{1}{4}, \\ \\ x_2^* = & \frac{4}{\tau - 4} &,\ 0< \tau \neq 4.\end{array} \right.
\end{eqnarray}

\begin{figure}[h]
\begin{center}
\begin{tikzpicture}[scale=0.9]
\draw [-triangle 45](0,0) node [left]{{\scriptsize -15}} -- (0,7);\draw (0,7) node [left]{{\footnotesize $x_2$}} -- (12,7) -- (12,0);\draw [-triangle 45] (0,0) -- (12,0) node [below]{{\footnotesize $x_1$}};
\draw[dashed] (0,3) node [left]{{\scriptsize 0}} -- (12,3);\draw[dashed] (6,0) node [below]{{\scriptsize 0}} -- (6,7);

\draw (0,1) node [left]{{\scriptsize -10}} -- (0.1,1);\draw (0,2) node [left]{{\scriptsize -5}} -- (0.1,2);\draw (0,4) node [left]{{\scriptsize 5}} -- (0.1,4);\draw (0,5) node [left]{{\scriptsize 10}} -- (0.1,5);\draw (0,6) node [left]{{\scriptsize 15}} -- (0.1,6);

\draw (0,0) node [below]{{\scriptsize -15}};
\draw (2,0.1) -- (2,0) node [below]{{\scriptsize -10}};\draw (4,0.1) -- (4,0) node [below]{{\scriptsize -5}};\draw (8,0.1) -- (8,0) node [below]{{\scriptsize 5}};\draw (10,0.1) -- (10,0) node [below]{{\scriptsize 10}};

\draw (11.9,1) -- (12,1);\draw (11.9,2) -- (12,2);\draw (11.9,4) -- (12,4);\draw (11.9,5) -- (12,5);\draw (11.9,6) -- (12,6);

\draw (2,6.9) -- (2,7) ;\draw (4,6.9) -- (4,7);\draw (8,6.9) -- (8,7);\draw (10,6.9) -- (10,7);


\draw[line width=1pt, dotted] (11.25,2.85) -- (11.6,2.85) node [below]{$A$};
\draw[line width=1pt] (11.2,2.85) .. controls (7.2,2.85) .. (6.65,2.92);


\draw[line width=1pt, dotted] (5.3,6.35) -- (5.3,6.7) node [left]{$B$};
\draw[line width=1pt] (5.4,3.7) .. controls (5.3,4) .. (5.3,6.3);

\draw[line width=1pt, dotted] (0.6,2.8) node [below]{$C$}-- (1,2.8);\draw[line width=1pt, dotted] (5.2,0.29) -- (5.2,0);
\draw[line width=1pt] (5.2,0.3) .. controls (5.1,2.8) .. (1,2.8);


\draw[] (5,5.5) node [left]{\scriptsize $\tau > 4$};
\draw[] (4,2) node [left]{\scriptsize $\frac{1}{4}< \tau < 4$};
\draw[] (10.5,2.3) node [left]{\scriptsize $0< \tau < \frac{1}{4}$};
\end{tikzpicture}
\caption{Locally Pareto-efficient optimal solutions for the Example~\ref{exemp:Exemplo1}}\label{fig:ramas}
\end{center}
\end{figure}
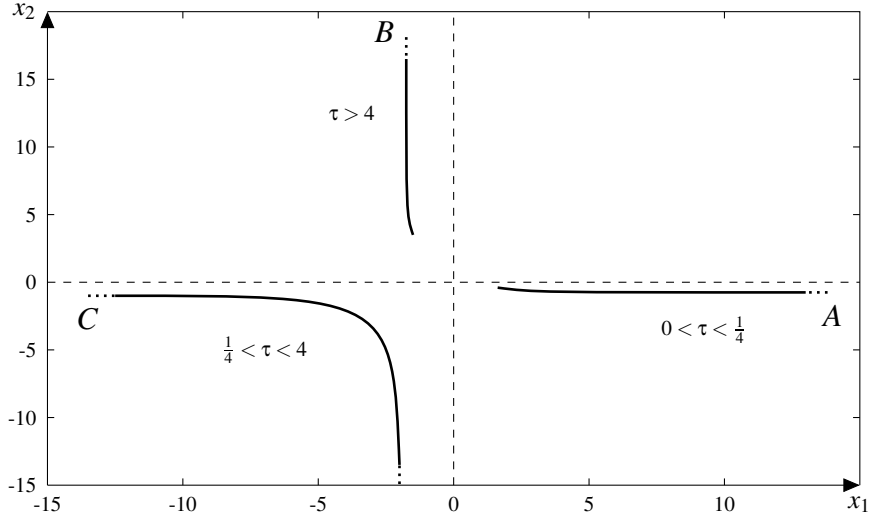

\begin{figure}[h]
\begin{center}
\begin{tikzpicture}[scale=0.9]
\draw [-triangle 45] (0,2) node [left]{{\scriptsize -50}} -- (0,9);
\draw (0,9) node [left]{{\footnotesize $f_2$}} -- (12,9) -- (12,2) node [below]{{\footnotesize $f_1$}}; 
\draw [-triangle 45](0,2) -- (12,2);
\draw[dashed] (0,3) node [left]{{\scriptsize 0}} -- (12,3);
\draw[dashed] (6,2) node [below]{{\scriptsize 0}} -- (6,9);

\draw (0,4) node [left]{{\scriptsize 50}} -- (0.1,4);\draw (0,5) node [left]{{\scriptsize 100}} -- (0.1,5);\draw (0,6) node [left]{{\scriptsize 150}} -- (0.1,6);\draw (0,7) node [left]{{\scriptsize 200}} -- (0.1,7);\draw (0,8) node [left]{{\scriptsize 250}} -- (0.1,8);

\draw (0,2) node [below]{{\scriptsize -150}};
\draw (2,2.1) -- (2,2) node [below]{{\scriptsize -100}};
\draw (4,2.1) -- (4,2) node [below]{{\scriptsize -50}};\draw (8,2.1) -- (8,2) node [below]{{\scriptsize 50}};
\draw (10,2.1) -- (10,2) node [below]{{\scriptsize 100}};

\draw (11.9,2) -- (12,2);\draw (11.9,4) -- (12,4);\draw (11.9,5) -- (12,5);\draw (11.9,6) -- (12,6);\draw (11.9,7) -- (12,7);\draw (11.9,8) -- (12,8);\draw (11.9,9) -- (12,9);

\draw (2,8.9) -- (2,9) ;\draw (4,8.9) -- (4,9);\draw (8,8.9) -- (8,9);\draw (10,8.9) -- (10,9);

\draw[] (10.7,2.6) node [below]{$f(A)$};
\draw[line width=1pt] (12,2.5) .. controls (7.3,2.9) .. (6.65,2.92);


\draw[] (5.6,7) node [left]{$f(B)$};
\draw[line width=1pt] (5.94,3.08) .. controls (5.5,4.5) .. (3.62,9);

\draw[] (1.7,8.5) node [below]{$f(C)$};
\draw[line width=1pt] (1.92,9) -- (2.97,7);
\draw[line width=1pt] (10,2.1) .. controls (5.6,2.49) .. (2.97,7);

\end{tikzpicture}
\caption{Pareto-optimal curve for the Example~\ref{exemp:Exemplo1}}\label{fig:Exemplo1Pareto}
\end{center}
\end{figure}
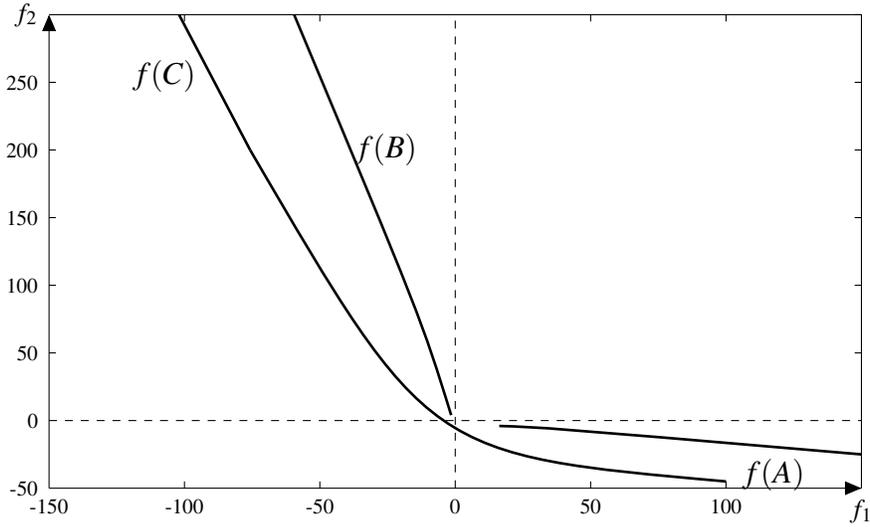

The points $(x_1^*,x_2^*)$ satisfying~(\ref{eq:PontosEficientesExemplo1}) are represented in Figure~\ref{fig:ramas}. The locally 
Pareto-efficient optimal solutions to this problem are grouped into three disconnected branches: $A$, $B$ and $C$. In addition, it 
is possible to verify that only the solutions set generated by inequality $\frac{1}{4} < \tau < 4$ are globally Pareto-efficient 
optimal solutions (branch $C$ plotted in the graph). We accept this by looking Figure~\ref{fig:Exemplo1Pareto}. It illustrates the graph in 
the plane representing the Pareto-optimal curve, that is, the branches images $A$, $B$ and $C$, by functions $f_1$ and $f_2$ plotted in 
the coordinates $(f_1(x_1^*, x_2^*), f_2(x_1^*, x_2 ^*))$. $f(A)=(f_1(A),f_2 (A))$, $f(B)=(f_1 (B),f_2(B))$ and $f(C)=(f_1(C),f_2(C))$ 
represent the branches images $A$, $B$ and $C$, respectively. We see in Figure~\ref{fig:Exemplo1Pareto} that the locally Pareto-efficient 
optimal solutions set belonging to the branch $C$ also are globally Pareto-efficient optimal solutions, that is, globally non-dominated.

We have just seen a simple example which the globally Pareto-efficient optimal solutions set is strictly contained in the locally 
Pareto-efficient optimal solutions set, but investigate the locally or globally efficient optimal solutions are not a trivial task. In 
this paper, we present new ways that can facilitate this analysis.

This paper is organized as follows. We start by defining some notations and basic properties. In Section~\ref{cap:Capitulo3}, the new concept of radius of 
efficiency and the relationships among associated problem are presented. In Section~\ref{cap:Capitulo4}, using this concept some necessarily and sufficient Pareto optimality conditions 
are established. Finally, comments and concluding remarks are presented In Section~\ref{cap:Conclusions}.

\section{Quadratic fractional problem and literature review}\label{cap:Capitulo2}

In this paper, we deal with a particular case of (VOP), which each objective function consists 
of a ratio of two quadratic functions. Without generalized convexity assumptions in the objective functions, we show how to calculate the largest radius of the spherical region
centered on a local Pareto-efficient solution in which this solution is optimal. Let us consider the following vector quadratic fractional optimization problem:

\begin{center}
\begin{tabular}{lll}
\textbf{(VQFP)}&\quad \mbox{\textbf{Minimize}}\quad &$\frac{f(x)}{g(x)}=\left( \frac{f_1(x)}{g_1(x)},\ldots,\frac{f_m(x)}{g_m(x)} \right)$\\
&\quad \mbox{subject to}& $h_j(x)\leqq 0,\quad j\in J$,\\
&                 & $x\in \Omega\subseteq \mathbf{R}^n$,
\end{tabular}
\end{center}
where $f_i$ and $g_i$, $i\in I$, are quadratic functions of $n$ real variables. In addition, we assume that $g_i(x)>0$, $i\in I$, 
for any $x$ in nonempty $\Omega$. The fields of real functions $h_j's$ contains $\Omega$. We denote by $S$ the \textit{feasible set} 
that is the intersection of $\Omega$ with the set of points $x$ in which $h_j(x) \leqq 0$, $j\in J$. We choose the functions $g_i's$ 
that preserves the signal. We denote the unconstrained (VQFP) by (VQFP$'$).

The fractional optimization problems arise frequently in the decision making applications, including science management, 
portfolio selection, cutting and stock, game theory, in the optimization of the ratio performance/cost, or profit/ investment, or cost/time and so on. 

There are many contributions dealing with the scalar (single-objective) fractional optimization problem (FP) and vector fractional optimization 
problem (VFP). In most of them, using convexity or their generalizations, optimality conditions in the KKT sense and the main duality theorems for optimal points are obtained. With a parametric approach, which transforms the original problem in a simpler associated problem, Dinkelbach~\cite{Art:Dinkelbach}, Jagannathan~\cite{Art:Jagannathan} 
and Antczak~\cite{Art:Antczak} establish optimality conditions, presents algorithms and apply their approaches in a example (FP) consisting 
of quadratic functions. For example, Dinkelbach~\cite{Art:Dinkelbach} presents some theoretical results that relate the two problems, and 
proposes an algorithm that converges to the minimum of problem (FP) to perform a sequence of operations in associated parametric problem. Using 
some known generalized convexity, Antczak~\cite{Art:Antczak}, Khan and Hanson~\cite{Art:KhanAndHanson}, Reddy and Mukherjee~\cite{Art:ReddyAndMukherjee}, 
Jeyakumar~\cite{Art:Jeyakumar}, Liang et al.~\cite{Art:LiangAndHuangAndPardalos} establish optimality conditions and theorems that relate the pair primal-dual of problem (FP). In Craven~\cite{Art:Craven1981,Livro:Craven1988} and Weir~\cite{Art:Weir}, other results for the scalar optimization (FP) can be found.

Further, Liang et al.~\cite{Art:LiangAndHuangAndPardalosMultiobjetivo} extended their approach to the vector optimization case (VFP) considering the type duals of 
Mond-Weir~\cite{Art:MondWeirDual}, Schaible~\cite{Art:SchaibleDual1,Art:SchaibleDual2} and Bector~\cite{Art:BectorDual}. Considering the parametric approach of Dinkelbach~\cite{Art:Dinkelbach}, 
Jagannathan~\cite{Art:Jagannathan}, Bector et al.~\cite{Art:BectorAndChandraAndSingh} and two classes of generalized convexity, Osuna-G\'omez et al.~\cite{Art:Osuna} 
establish weak Pareto-efficient optimality conditions and the main duality theorems for the differenciable vector optimization case (VFP). Santos et al.~\cite{Art:Lucelina} deepened these results to the more 
general non-differenciable case (VFP). Jeyakumar and Mond~\cite{Art:JeyakumarAndMond} use generalized convexity to study the problem (VFP) and with the parametric 
approach, Singh and Hanson~\cite{Art:SinghAndHanson} extended the results obtained by Geoffrion~\cite{Art:Geoffrion}. 

Few studies are found involving quadratic functions at both the numerator and denominator in the ratio objective function. Most of them involve the mixing 
of linear and quadratic functions. The closest approaches of the scalar quadratic fracional optimization problem (QFP) are considered by Crouzeix et al.~\cite{CrouzeixAndFerlandAndSchaible}, Schaible and Shi~\cite{Art:SchaibleAndShiQuadratic}, 
Gotoh and Konno~\cite{Art:GotohAndKonno}, Lo and MacKinlay~\cite{Art:LoAndMacKinlay} and Cambini et al.~\cite{CambiniAndCrouzeixAndMartein}. On the other hand, Benson~\cite{Art:BensonQuadratic} considered a pure (QFP) consisting of the convex function, develop some theoretical properties and optimality conditions, 
he presents an algorithm and its convergence properties. 

The closest approaches of the vector optimization case (VQFP) are considered by Beato et al.~\cite{Art:BeatoAndRuizAndLuqueAndBlanquero,Art:BeatoAndInfanteAndRuiz}, 
Ar\'evalo and Zapata~\cite{Art:ArevaloAndZapata}, Konno and Inori~\cite{KonnoAndInori}, Rhode and Weber~\cite{Livro:RhodeAndWeber}, Kornbluth and Steuer~\cite{Art:KornbluthAndSteuer}, 
Korhoen and Yu~\cite{Art:KorhonenAndYu,Art:KorhonenAndYu2}. Using an iterative computational test, 
Beato et al.~\cite{Art:BeatoAndRuizAndLuqueAndBlanquero,Art:BeatoAndInfanteAndRuiz} characterize the Pareto-efficient optimal point for the 
problem (VQFP), consisting of a linear and quadratic functions, and using the function linearization approach of Bector et al.~\cite{Art:BectorAndChandraAndSingh}, some theoretical results are obtained. Ar\'evalo and Zapata~\cite{Art:ArevaloAndZapata}, Konno and Inori~\cite{KonnoAndInori}, Rhode and Weber~\cite{Livro:RhodeAndWeber} analyze the 
portfolio selection problem. Kornbluth and Steuer~\cite{Art:KornbluthAndSteuer} use an adapted Simplex method in the problem (VFP) consisting of linear functions. 
Korhonen and Yu~\cite{Art:KorhonenAndYu,Art:KorhonenAndYu2} propose an iterative computational method for solving the problem (VQFP), consisting of the linear and 
quadratic functions, based on search directions and weighted sums.

The approach taken in this work is different from the previous ones. We believe that the approach presented here facilitates the resolution of the problem (VQFP). 
The main aspect of this contribution is the development of necessary and sufficient Pareto-efficient optimality conditions for a particular 
vector optimization problems based on the calculation of largest radius of the spherical region centered on a local Pareto-efficient solution in which 
this solution is optimal. In this process we may conclude that the solution is also globally Pareto-efficient optimal. These conditions might be useful to 
determine termination criteria in the development of algorithms, including more general problems in which quadratic approximations are used locally.

\section{Radius of efficiency}\label{cap:Capitulo3}

We introduce the new concept \textit{radius of efficiency} for (VOP) and present some Pareto-efficient optimality conditions for (VQFP$'$) 
from this alternative approach, and then we extend the results to the constrained case (VQFP). 

\subsection{Radius of efficiency in the (VOP)}

From certain particular properties verified in some search directions of the objective functions, we want to detect when a local Pareto-efficient 
optimal solution is equivalent to a global Pareto-efficient optimal solution. With this approach is possible to identify a spherical 
region of feasible points where a local Pareto-efficient optimal solution is not dominated.

We say that a point $x^*\in S$ is \textit{$\lambda$-efficient} or has \textit{radius of efficiency $\lambda$} in (VOP) if 
$x^*\in$ \textit{Leff}$(VOP)$ and there does not exist another point $x'\in B(x^*,\lambda)\cap S$ that dominates $x^*$.

Note an important difference between the locally Pareto-efficient optimal solution definition and the radius of efficiency definition of a locally 
Pareto-efficient optimal solution. In the first case, from a theoretical point of view we know that always there is an arbitrary neighbourhood 
$N(x^*)\subseteq S$ of $x^*$, where $x^*$ is not dominated. Naturally, this neighbourhood always can be regarded as a ball of arbitrary radius 
in $\R^n$. In the second case, is possible to calculate the largest radius $\lambda^*>0$ of the spherical region $B(x^*,\lambda^*)$ 
such that $x^*$ is not dominated in $B(x^*,\lambda^*)\cap S$, or is possible to conclude that $x^*$ is not dominated everywhere 
the feasible set $S$.

Some important reasons for the use of the radius of efficiency in the problems (VOP) can be indicate. We always can consider a local solution in 
a fixed well-know spherical region instead in an arbitrary spherical region. We can determine a subset compact where there does not exist 
another points that dominate a specific solution, what is useful in the resolution of (VOP), because if the decision maker knows the radius of 
efficiency, he can estimates the cost to try to find a new solution, and also choose a more suitable search procedure. Auxiliary problems induced 
by the concept of radius of efficiency can be used to conclude the global Pareto-efficiency.

Naturally, if $x^*$ is $\lambda$-efficient, then it is $\beta$-efficient, $\forall\ \beta < \lambda$. Similarly we say that $x^*$ is $\infty$-efficient if it is efficient in $S$.

\subsection{Radius of efficiency in the (VQFP)}\label{sec:raioEmPMFQ}

We calculate the radius of efficiency of a solution $x^*$ and some results from the calculation are presented. By hypothesis, we assume 
that $x^*\in$ \textit{Leff}$(VQFP)$. To determine the radius of efficiency $\lambda>0$ for this solution $x^*$, we must 
ask: \textit{What is the smallest value of $\lambda>0$ such that among every unitary direction $d$ in which $\frac{f(x^*+\lambda d)}{g(x^*+\lambda d)} \leq \frac{f(x^*)}{g(x^*)}$ is valid?} 
The answer to this question provides the maximum radius of efficiency of $x^*$. 

The next theorems allow us to characterise when a locally Pareto-efficient optimal solution is equivalent to a Pareto-efficient optimal solution. 
We will use these theorems to identify the maximum radius of efficiency and analyse the dominance of a locally Pareto-efficient optimal solution 
in the feasible set.

Similarly to Dinkelbach~\cite{Art:Dinkelbach} and Jagannathan~\cite{Art:Jagannathan}, which transform the fractional optimization in a new problem, 
we consider the following problem associated with the (VQFP).

{\footnotesize \begin{eqnarray*}
\textbf{(VQFP)}_{x^*} &\ \mbox{\textbf{Minimize}} &\ \ f(x)-\frac{f(x^*)}{g(x^*)}g(x)=\left( f_1(x)-\frac{f_1(x^*)}{g_1(x^*)}g_1(x),\ldots,f_m(x)-\frac{f_m(x^*)}{g_m(x^*)}g_m(x) \right)\\
&\ \mbox{subject to}  &\ \ h_j(x)\leqq 0\quad j\in J,\\
&                 & \ \ x\in \Omega\subseteq \mathbf{R}^n,
\end{eqnarray*}}were $x^*\in S$ and $f_i$, $g_i$, $i\in I$, $h_j$, $j\in J$, are the same functions defined in (VQFP).\\

The (VQFP)$_{x^*}$ is similarly introduced by Osuna-G\'omez et al., which derive optimality conditions and duality results for the weakly 
Pareto-efficient optimal solutions. The following theorem and its proof is an approach equivalent to Lemma 1.1 presented in~\cite{Art:Osuna}, 
but we consider the Pareto-efficient optimal solutions.

{\tteo{$x^*\in$ \textit{Leff}$(VQFP)$ if and only if $x^*\in$ \textit{Leff}$(VQFP)_{x^*}$. In addition, $x^*$ is locally Pareto-efficient optimal 
solution for $(VQFP)$ in $N(x^*)$ if and only if is locally Pareto-efficient optimal solution for $(VQFP)_{x^*}$ in $N(x^*)$.}\label{teo:ProblemaAssociado}}\\

{\prova{$(\Rightarrow)$ Let $N(x^*)\subseteq S$ and $x^*$ be locally Pareto-efficient optimal solution for $(VQFP)$ in $N(x^*)$. Suppose 
that $x^*\notin$ \textit{Leff}$(VQFP)_{x^*}$, then there exists another point $x'\in N(x^*)$ satisfying
\begin{eqnarray*}
f(x')-\frac{f(x^*)}{g(x^*)}g(x')\leq f(x^*)-\frac{f(x^*)}{g(x^*)}g(x^*)=0\quad \Longrightarrow\quad \frac{f(x')}{g(x')} \leq \frac{f(x^*)}{g(x^*)}.
\end{eqnarray*}
Which contradicts $x^*\in$ \textit{Leff}$(VQFP)$ in $N(x^*)$, and therefore $x^*\in$ \textit{Leff}$(VQFP)_{x^*}$ in $N(x^*)$.

$(\Leftarrow)$ Similarly, let $N(x^*)\subseteq S$ and $x^*$ be locally Pareto-efficient optimal solution for $(VQFP)_{x^*}$ in $N(x^*)$. 
Suppose that $x^*\notin$ \textit{Leff}$(VQFP)$, then there exists another point $x'\in N(x^*)$, satisfying
\begin{eqnarray*}
\frac{f(x')}{g(x')} \leq \frac{f(x^*)}{g(x^*)}\quad \Longrightarrow\quad f(x')-\frac{f(x^*)}{g(x^*)}g(x')\leq 0= f(x^*)-\frac{f(x^*)}{g(x^*)}g(x^*).
\end{eqnarray*}
Which contradicts $x^*\in$ \textit{Leff}$(VQFP)_{x^*}$ in $N(x^*)$, and therefore $x^*\in$ \textit{Leff}$(VQFP)$ in $N(x^*)$. }}\\

Note that on associated problems, for each $x^*\in$ \textit{Leff}$(VQFP)$ we can only ensure that \textit{Leff}$(VQFP)\subset$ \textit{Leff}$(VQFP)_{x^*}$. 

The same arbitrary neighborhood $N(x^*)$ equivalent to a ball of arbitrary radius in $\R^n$ centered at $x^*$ appears in both problems (VQFP) and 
(VQFP)$_{x^*}$, and we want to compute a neighborhood $B(x^*,\lambda^*)$ of maximum radius $\lambda^*$, such that $x^*\in$ \textit{Leff}$(VQFP)$ 
in $B(x^*,\lambda^*)$. Then, calculate the radius of efficiency $\lambda^*$ of the solution $x^*$ in the (VQFP) is equivalent to calculate the 
radius of efficiency $\lambda^*$ of the solution $x^*$ in the (VQFP)$_{x^*}$, and so we can choose between the two problems one whose 
calculation is easier.

\section{Optimality Conditions from Radius of efficiency}\label{cap:Capitulo4}

First the optimality conditions are established for the unconstrained problem and then we extend the results for the feasible set defined by linear 
inequalities.

\subsection{Radius of efficiency in the unconstrained case}\label{sec:raioEmQFMP}

For each $i\in I$ and all $x\in \mathbf{R}^n$ we consider the objective functions defined as $f_i(x)=x^T A_i x + a_i^T x + \bar{a}_i$ and $g_i(x)= x^T B_i x + b_i^T x + \bar{b}_i,$
where $A_i$, $B_i\in \mathbf{R}^{n\times n}$, $A_i$ symmetric, $B_i$ symmetric and positive semidefinite, $a_i$, $b_i\in \mathbf{R}^n$ and $\bar{a}_i$, $\bar{b}_i\in \mathbf{R}$, with
$\bar{b}_i > -( {w}^{i^T} B_i w^i + b_i^T w^i)$, where $w^i$ is the solution of
the system $2B_i x + b_i=0$, that is, $w^i$ is the point where the function $x^T B_i x + b_i^T x$ reaches its minimum and this ensures 
that $g_i(x)>0$, $\forall x\in \mathbf{R}^n$. We cannot consider cases where $2 B_i x + b_i= 0$ has no solution. 
Similarly, we denote by (VQFP$'$)$_{x^*}$ the unconstrained (VQFP)$_{x^*}$. We recall that $\nabla^2 f(x)$ denote the Hessian matrix of the function 
$f:\mathbf{R}^n\to \mathbf{R}$ at point $x$.

Further, we define some sets and parameters that are used throughout this work. Given $x^*\in \R^n$, we define the following quadratic function
\begin{eqnarray}\label{eq:DefineFuncoesP}
p_i(x)=x^T\left( A_i  - \frac{f_i(x^*)}{g_i(x^*)}  B_i \right)x + \left(a_i^T - \frac{f_i(x^*)}{g_i(x^*)} b_i^T\right) x,\quad i\in I,
\end{eqnarray}
and given an unitary direction $d$, we define

$$X_0\label{eq:ConjuntosX012}=\{i\in I\ |\ \left[d^T \nabla^2 p_i(x^*) d,\ \nabla p_i(x^*)^T d\right]^T \geq \left[0,0\right]^T\},$$

$$X_1=\{i\in I\ |\ d^T \nabla^2 p_i(x^*) d >0\ \ \mbox{and}\ \ \nabla p_i(x^*)^T d < 0\},$$

$$X_2=\{i\in I\ |\ d^T \nabla^2 p_i(x^*) d <0\ \ \mbox{and}\ \ \nabla p_i(x^*)^T d > 0\},$$

$$\lambda_{2}^d = \max\limits_{i \in X_2} \left\{\ \frac{-2\nabla p_i(x^*)^T d}{d^T \nabla^2 p_i(x^*) d}\ \right\}, \qquad \lambda_1^d = \left\{\begin{tabular}{ll}$\min\limits_{i \in X_1} \left\{\displaystyle \frac{-2\nabla p_i(x^*)^T d}{d^T \nabla^2 p_i(x^*) d} \right\}$,& \mbox{if}\ $X_1\neq \emptyset$ \\
$+\infty$, & \mbox{if}\ $X_1 = \emptyset,$ \end{tabular}\right.$$

$$\Lambda_{2}^d=[\lambda_{2}^d,\ \infty),\qquad\qquad \Lambda_1^d= \left\{\begin{tabular}{ll}$(0,\ \lambda_1^d ]$, & \mbox{se}\ $X_1\neq\emptyset$ \\ $(0,\ \lambda_1^d )$, & \mbox{se}\ $X_1=\emptyset$, \end{tabular}\right.\qquad \quad \Lambda^d = \Lambda_{2}^d \cap \Lambda_1^d.$$

The functions $p_i$, $i\in I$, defined in~(\ref{eq:DefineFuncoesP}) are the same objective functions of the (VQFP)$_{x^*}$ unless the constants 
term $ \bar{a}_i - \frac{f_i(x^*)}{g_i(x^*)} \bar{b}_i$, and so $ p_i(x^*)=-[\bar{a}_i - \frac{f_i(x^*)}{g_i(x^*)} \bar{b}_i]$. We denote by $|A|$ the 
number of elements in the set $A$.

For a better understanding of sets and parameters above, consider the Taylor expansion around zero of the function 
$r_i(\lambda)=p_i(x^* + \lambda d)$, $i\in I$, $\lambda\in \R$. Since $r_i$ is a quadratic function, we obtain \begin{displaymath} r_i(\lambda) = p_i(x^*) + \lambda \nabla p_i(x^*)^T d + \frac{\lambda^2}{2} d^T \nabla^2 p_i(x^*) d. \end{displaymath} 
In other words, $r_i$ is a real function of one variable, whose graph is a parable, and has constant term $p_i(x^*)$, linear term 
$\nabla p_i(x^*)^T d$ and quadratic term $\frac{1}{2} d^T \nabla^2 p_i(x^*) d$. Thus, it becomes easier to interpret the sets $X_0$, $X_1$ and $X_2$. 
Given an unitary direction $d$, the set $X_1$ is formed by the indices $i\in I$ whose $r_i$ functions decrease to values close $\lambda=0$ and are 
convex, while the set $X_2$ is formed by the indices $i\in I$ whose $r_i$ functions grow to values close $\lambda=0$ and are concave. 

\begin{figure}[h]
\begin{center}
\begin{tikzpicture}

\draw [->](1,0.5) -- (1,5) node [above]{\small $r_i(\lambda)$};
\draw [->](.5,1) -- (7,1) node [right]{\small $\lambda$};

\draw[color=black, line width=.7pt] (1,3) sin (3,1.5) cos (6.6,5);
\draw [dashed](5.2,3) -- (5.2,.7);
\draw (5,1) node[above] {{\scriptsize \textcolor{black}{$\lambda_1^{d}$}}};
\draw [color=black,<-|, line width=1.2pt] (1,1) -- (5.2,1);
\draw (1.4,1) node[above] {{\scriptsize \textcolor{black}{$\Lambda_1^{d}$}}};
\draw [color=black,<->, line width=1.2pt] (3.6,.85) -- (5.2,.85);\draw (4.4,.85) node[below] {{\scriptsize \textcolor{black}{$\Lambda^{d}$}}};

\draw [dashed](1,3) -- (7,3) (0.5,3) node {{\scriptsize $p_{i} (x^*)$}};
\draw [line width=.7pt] (1,3) parabola[bend pos=0.5] bend +(0,2.5) +(3,-1.5);
\draw [|->, line width=1.2pt](3.58,3) -- (7,3);
\draw (6.8,3.05) node[above] {{\scriptsize \textcolor{black}{$\Lambda_2^{d}$}}};
\draw [dashed](3.6,3) -- (3.6,0.7);
\draw (3.75,3.05) node[above] {{\scriptsize $\lambda_{2}^{d}$}};


\draw [<-](2.9,4.6) -- (3.5,5);
\draw (3.5,5.2) node {$i\in X_2$};

\draw [color=black,->](5.55,5) -- (6.2,4.6);
\draw (5.5,5.2) node [color=black]{$i\in X_1$};

\end{tikzpicture}
\caption{Graphs of $r_i$ in which index $i$ belongs to the set $X_1$ or $X_2$}\label{fig:ElementosEmX1eX2}
\end{center}
\end{figure}
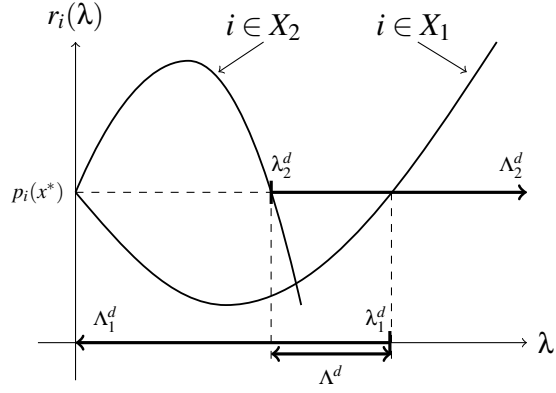

\begin{figure}[h]
\begin{center}
\begin{tikzpicture}

\draw [->](1,0.5) -- (1,5) node [above]{\small $r_i(\lambda)$};
\draw [->](.5,1) -- (7,1) node [right]{\small $\lambda$};


\draw [dashed](1,3) -- (7,3) (0.5,3) node {{\scriptsize $p_{i} (x^*)$}};
\draw [line width=.8pt] (1,3) parabola[bend pos=0.5] bend +(0,2) +(5,-1.5);

\draw[color=black, line width=.8pt] (1,3) -- (4.5,5.6);
\draw[color=black, line width=.8pt] (1,3) cos (4.2,5);

\draw [<-](5.3,3.1) -- (6.3,3.85);
\draw (6.3,4.25) node {$\bar{\lambda}=\frac{-2\nabla p_i(x^*)^T d}{d^T \nabla^2 p_i(x^*) d}$};

\draw [<-, color=black](4.2,4.85) -- (5.1,5.2);
\draw [color=black] (5.7,5.2) node {$i\in X_0$};

\draw [<-, color=black](4.5,5.5) -- (5.1,5.2);

\draw [->](4.65,1.8) -- (5.55,2);
\draw (4,1.8) node {$i\in X_2$};

\end{tikzpicture}
\caption{Graphs of $r_i$ in which index $i$ belongs to the set $X_0$ or $X_2$, and a root $\bar{\lambda}$ of the equation $r_i(\lambda)-p_i(x^*)=0$}\label{fig:ElementosEmX0eX2eRaiz}
\end{center}
\end{figure}
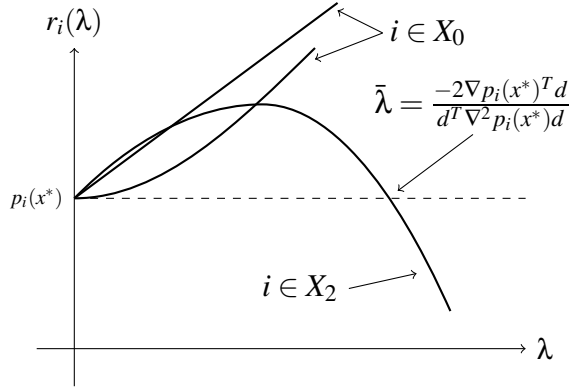

Figure~\ref{fig:ElementosEmX1eX2} shows typical elements of the sets $X_1$ and $X_2$, in which the graphs represent two functions $r_i$, and are 
plotted the points of function $r_i(\lambda)=p_i(x^* + \lambda d)$ in the coordinates $(\lambda,r_i(\lambda))$ for $\lambda>0$, which are two 
parables: a convex indicating that the function index $r_i$ belongs to the set $X_1$ (where, $d^T \nabla^2 p_i(x^*) d>0$ and 
$\nabla p_i(x^*)^T d<0$); and a concave indicating that the function index $r_i$ belongs to the set $X_2$ 
(where, $d^T \nabla^2 p_i(x^*) d<0$ and $\nabla p_i(x^*)^T d>0$).

The set $X_0$ is formed by the indices $i\in I$ where each quadratic function $r_i$ has non-negative linear term and positive quadratic term 
or has positive linear term and non-negative quadratic term. 

Figure~\ref{fig:ElementosEmX0eX2eRaiz} shows typical elements of the set $X_0$ and $X_2$, in which the graphs represent three functions $r_i$. In 
relation the function index $r_i$ belongs to the set $X_0$, two examples are plotted: a convex parable (where, $d^T \nabla^2 p_i(x^*) d>0$ 
and $\nabla p_i(x^*)^T d \geqq 0$); and a growing straight line (where, $d^T \nabla^2 p_i(x^*) d=0$ and $\nabla p_i(x^*)^T d > 0$).

The parameter $\lambda_2^d$ is defined as the maximum of the positive roots of the equations $r_i(\lambda)-p_i(x^*)=0$ in which indices $i$ belong 
to $X_2$. Similarly, if $X_1\neq \emptyset$, the parameter $\lambda_1^d$ is defined as the minimum of the positive roots of the equations 
$r_i(\lambda)-p_i(x^*)=0$ in which indices $i$ belong to $X_1$. If $X_1=\emptyset$, we define $\lambda_1^d\equiv +\infty$. 

For each index $i$, another root of the equation $r_i(\lambda)-p_i(x^*)=0$ is $\lambda=0$, but we are interested in the values $\lambda>0$. 
Figure~\ref{fig:ElementosEmX1eX2} shows these parameters, the parameter $\lambda_2^d$ and $\lambda_1^d$ are the crossover point between 
the curve $r_i(\lambda)$ and the dotted straight line $p_i(x^*)$ for $i\in X_2$ and $i\in X_1$, respectively. Figure~\ref{fig:ElementosEmX0eX2eRaiz} shows a positive root of the equation $r_i(\lambda)-p_i(x^*)$ for $i\in X_2$.
\begin{figure}[h]
\begin{center}
\begin{tikzpicture}

\draw [->](1,0.5) -- (1,5) node [above]{\small $r_i(\lambda)$};
\draw [->](.5,1) -- (7,1) node [right]{\small $\lambda$};

\draw[color=black, line width=.8pt] (1,3) sin (2.2,1.5) cos (4.5,5);

\draw [color=black, line width=.8pt, dashed](1,3) -- (7,3);\draw (0.5,3) node {{\scriptsize $p_{i} (x^*)$}};

\draw[color=black, line width=.8pt] (1,3) -- (4.5,0.5);
\draw[color=black, line width=.8pt] (1,3) cos (4,1.1);

\draw [<-, color=black](3.65,2.92) -- (4.56,2.15);
\draw [color=black] (5.9,2.1) node {$\bar{\lambda}=\frac{-2\nabla p_i(x^*)^T d}{d^T \nabla^2 p_i(x^*) d}$};

\draw [<-, color=black](4.5,4.4) -- (5.4,4.5);
\draw [color=black] (6,4.5) node {$i\in X_1$};

\draw [<-, color=black](4.1,1.16) -- (4.55,1.2);
\draw [<-, color=black](4.3,0.75) -- (4.55,1.2);
\draw [color=black] (5.9,1.35) node {\scriptsize items 2(c), 3(b) and 3(c)};

\draw [<-, color=black](5.4,3.09) -- (5.8,3.2);
\draw [color=black] (6.4,3.3) node {\scriptsize item 2(b)};

\end{tikzpicture}
\caption{Graphs of $r_i$ in which index $i$ belongs to the set $X_1$, and a root $\bar{\lambda}$ of the equation $r_i(\lambda)-p_i(x^*)=0$}\label{fig:ElementosEmX1eRaiz}
\end{center}
\end{figure}
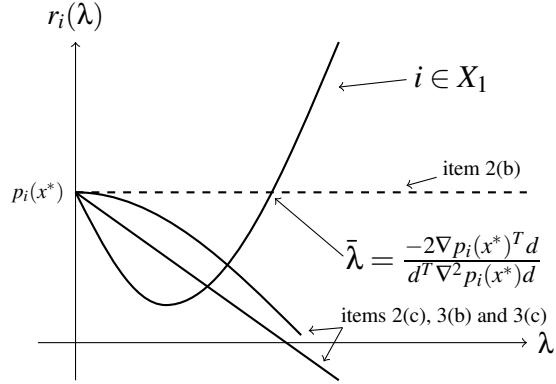

Given an unitary direction $d$, the sets $\Lambda_{2}^d$, $\Lambda_{1}^d$ and $\Lambda^d$ are real intervals contained in $\R_{+}\setminus\{0\}$, 
where $\lambda_{2}^d$ is the extreme left value of the interval $\Lambda_{2}^d$ and $\lambda_{1}^d$ is the extreme right value of the 
interval $\Lambda_{1}^d$. When occur $\lambda_{2}^d\leqq \lambda_{1}^d$, they are the end points of the interval $\Lambda^d$ and if occur $\lambda_{2}^d > \lambda_{1}^d$, we 
obtain $\Lambda^d=\emptyset$. However, in some particular cases, for example when occur $|X_1|=|X_2|= 1$ and $\lambda_{2}^d = \lambda_{1}^d$, 
we obtain $\Lambda^d=\{\lambda_{2}^d\}$, that is $\Lambda^d\neq \emptyset$, but the most suitable in this case is set $\Lambda^d=\emptyset$. This 
case and others like it are important and are explained in detail in due course.

\begin{figure}[h]
\begin{center}
\begin{tikzpicture}

\draw [->](1,0.5) -- (1,5) node [above]{\small $r_i(\lambda)$};
\draw [->](.5,1) -- (7,1) node [right]{\small $\lambda$};

\draw [dashed](.9,2) -- (7,2) (0.5,2) node {{\tiny $p_{4} (x^*)$}}; \draw [line width=.7pt] (1,2) parabola[bend pos=0.5] bend +(0,4) +(4,-0.5);\draw [->, line width=1.2pt](4.9,2) -- (5.3,2);

\draw [dashed](.9,2.5) -- (7,2.5) (0.5,2.5) node {{\tiny $ \textcolor {black}{p_{3} (x^*)}$}};\draw[color=black, line width=.7pt] (1,2.5) sin (2.2,1.5) cos (4.5,5);\draw [color=black,<-, line width=1.2pt](3,2.5) -- (3.35,2.5);\draw [dashed](3.35,2.5) -- (3.35,.8) node[below] {{\scriptsize \textcolor{black}{$\lambda_1^{d}$}}};\draw [color=black,<-|, line width=1.2pt](2.6,1) -- (3.35,1);

\draw [dashed](.9,3) -- (7,3) (0.5,3) node {{\tiny $p_{2} (x^*)$}};\draw [line width=.7pt] (1,3) parabola[bend pos=0.5] bend +(0,2) +(5,-1.5);\draw [->, line width=1.2pt](5.2,3) -- (5.6,3);\draw [dashed](5.2,3) -- (5.2,0.8) node[below] {{\scriptsize $\lambda_{2}^{d}$}};\draw [|->, line width=1.2pt](5.2,1) -- (6,1);

\draw [dashed](.9,3.5) -- (7,3.5) (0.5,3.5) node {{\tiny $\textcolor {black}{p_{1} (x^*)}$}};\draw[color=black, line width=.7pt] (1,3.5) sin (2.5,2.5) cos (6,5.5);\draw [color=black,<-, line width=1.2pt](4,3.5) -- (4.4,3.5);

\draw (6.7,4.5) node {$d\notin L$};\draw (6.3,0.5) node {{\small $\Lambda^{d}=\emptyset $}};

\draw [<-, color=black](4.65,5) -- (5.3,5.4);
\draw [->, color=black](5.3,5.4) -- (5.55,5);
\draw [color=black] (5.3,5.6) node {\small $1,3\in X_1$};

\draw [<-](2,4.56) -- (2.57,4.7);
\draw [->](2.57,4.7) -- (2.5,4.16);
\draw [] (3,4.9) node {\small $2,4\in X_2$};

\end{tikzpicture}
\caption{A search direction by $x'$ that dominates $x^*\in$ \textit{Leff}$(VQFP')_{x^*}$, $X_0=\emptyset$, $X_1=\{1,3\}$, $X_2=\{2,4\}$ and $\Lambda^{d}=\emptyset$}\label{fig:Direcao1BuscaSemRestricoes}
\end{center}
\end{figure}
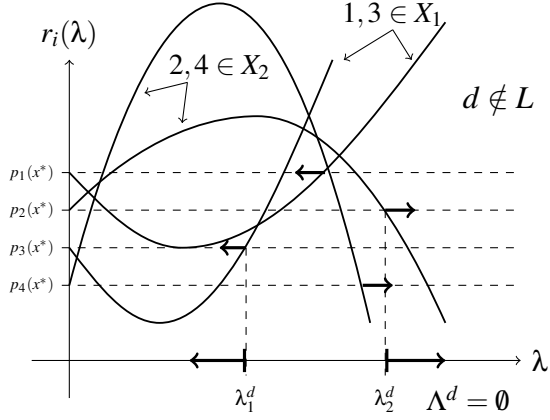

Figure~\ref{fig:ElementosEmX1eX2} illustrates examples of the intervals $\Lambda_{2}^d$, $\Lambda_{1}^d$ and $\Lambda^d$. In it, we represent over the 
dotted straight line $p_i(x^*)$ the interval $\Lambda_{2}^d=[\lambda_{2}^d,\infty)$, over 
the $\lambda$-axis the interval $\Lambda_{1}^d=(0,\lambda_{1}^d]$, and below the $\lambda$-axis the interval 
$\Lambda^d=\Lambda_{2}^d\cap \Lambda_{2}^d$. When $X_1=\emptyset$ we choose $\lambda_1^d=\infty$, because when $X_2\neq\emptyset$ we always obtain $\Lambda^d\neq \emptyset$. Note that in Figure~\ref{fig:ElementosEmX1eX2} 
we have $p_i(x')\leqq p_i(x^*)$ for all $\lambda\in \Lambda^d$ and $x'=x^* + \lambda d$. On the other hand, Figure~\ref{fig:Direcao1BuscaSemRestricoes} 
illustrates an example in which $\Lambda^d=\emptyset$.

Before the next result, we present some important details in order to understand the possibility of obtaining an unitary direction $d$ and a constant $\lambda>0$, in which 
$\frac{f(x^*+\lambda d)}{g(x^*+\lambda d)}\leq \frac{f(x^*)}{g(x^*)}$ is valid for a solution $x^*\in$\textit{Leff}$(VQFP')$, and what are their 
relations with the sets $X_0$, $X_1$ and $X_2$.

Consider the Taylor expansion around zero of the each function $\bar{r}_i(\lambda)=f_i(x^* + \lambda d)$ and 
$\tilde{r}_i(\lambda)=g_i(x^* + \lambda d)$, $i\in I$, in solution $x^*$ and along the direction $d$,
\begin{eqnarray*}
\bar{r}_i(\lambda)&=&f_i(x^*+\lambda d) = f_i(x^*) + \lambda \nabla f_i(x^*)^T d + \lambda^2 d^T A_i d,\\ 
\tilde{r}_i(\lambda)&=&g_i(x^*+\lambda d) = g_i(x^*) + \lambda \nabla g_i(x^*)^T d + \lambda^2 d^T B_i d.
\end{eqnarray*}
Performing some manipulations, we obtain 

{\small $$ \frac{f_i(x^*+\lambda d)}{g_i(x^*+\lambda d)} \leqq \frac{f_i(x^*)}{g_i(x^*)}\ \Longleftrightarrow \ \lambda d^T \left( \lambda A_i  - \lambda \frac{f_i(x^*)}{g_i(x^*)}  B_i \right)d \leqq \left(\lambda \frac{f_i(x^*)}{g_i(x^*)} \nabla g_i(x^*) - \lambda \nabla f_i(x^*)\right)^T d.$$}
From~(\ref{eq:DefineFuncoesP}) consider the function $p_i(x)$, then
\begin{eqnarray}
\lambda d^T \left( \lambda A_i  - \lambda \frac{f_i(x^*)}{g_i(x^*)}  B_i \right)d & \leqq & \left(\lambda \frac{f_i(x^*)}{g_i(x^*)} \nabla g_i(x^*) - \lambda \nabla f_i(x^*)\right)^T d\ \Longleftrightarrow 
\nonumber \\
\nonumber \\
& \Longleftrightarrow & \lambda\left( \nabla p_i(x^*)^T d + \frac{ \lambda }{2} d^T \nabla^2 p_i(x^*) d \right)  \leqq  0,
\label{eq:desig1}
\end{eqnarray}
where $ \frac{1}{2}\nabla^2 p_i(x^*) = A_i  - \frac{f_i(x^*)}{g_i(x^*)} B_i$ and $ \nabla p_i(x^*) = \nabla f_i(x^*) - \frac{f_i(x^*)}{g_i(x^*)} \nabla g_i(x^*).$ From~(\ref{eq:desig1}), 
we obtain 
\begin{eqnarray*}
p_i(x^* +\lambda d) - p_i(x^*) &=& \lambda\left( \nabla p_i(x^*)^T d + \frac{ \lambda }{2} d^T \nabla^2 p_i(x^*) d \right) \ \leqq \ 0 \ \Longleftrightarrow   \\
\nonumber \\
& \Longleftrightarrow & r_i(\lambda)\ =\ p_i(x^* +\lambda d) \ \leqq \ p_i(x^*).
\end{eqnarray*}
The following situations about the existence of positive solutions to $\lambda$ in 
inequalities~(\ref{eq:desig1}) can occur for each direction $d$ and for each $i\in I$.

\begin{enumerate}
\item { If $d^T \nabla^2 p_i(x^*) d >0$,
\begin{enumerate}
\item and if $\nabla p_i(x^*)^T d > 0$, does not exist $\lambda > 0$ satisfying~(\ref{eq:desig1}) and $i\in X_0$,
\item and if $\nabla p_i(x^*)^T d = 0$, does not exist $\lambda > 0$ satisfying~(\ref{eq:desig1}) and $i\in X_0$,  
\item and if $\nabla p_i(x^*)^T d < 0$, $ \lambda \in \left(0,\frac{-2\nabla p_i(x^*)^T d}{d^T \nabla^2 p_i(x^*) d}\right]$ 
satisfying~(\ref{eq:desig1}) and $i\in X_1$.
\end{enumerate}  
}

\item { If $d^T \nabla^2 p_i(x^*) d =0$,
\begin{enumerate}
\item and if $\nabla p_i(x^*)^T d > 0$, does not exist $\lambda > 0$ satisfying~(\ref{eq:desig1}) and $i\in X_0$,
\item and if $\nabla p_i(x^*)^T d = 0$, $\lambda \in (0,\infty)$ satisfying~(\ref{eq:desig1}),
\item and if $\nabla p_i(x^*)^T d < 0$, $\lambda \in (0,\infty)$ satisfying~(\ref{eq:desig1}).
\end{enumerate}
}

\item { If $d^T \nabla^2 p_i(x^*) d <0$,
\begin{enumerate}
\item and if $\nabla p_i(x^*)^T d > 0$, $ \lambda \in \left[\frac{-2\nabla p_i(x^*)^T d}{d^T \nabla^2 p_i(x^*) d},\infty\right)$ satisfying~(\ref{eq:desig1}) and $i\in X_2$,
\item and if $\nabla p_i(x^*)^T d = 0$, $\lambda \in (0,\infty)$ satisfying~(\ref{eq:desig1}),
\item and if $\nabla p_i(x^*)^T d < 0$, $\lambda \in (0,\infty)$ satisfying~(\ref{eq:desig1}).
\end{enumerate}
}
\end{enumerate}

The above items inform whenever exists $ \lambda > 0$ satisfying the inequality~(\ref{eq:desig1}) and in which interval it 
exists. They also inform the behavior of each function $r_i$ along the direction $d$, e.g. items 1.(a), 1.(b) and 2.(a) represent 
the functions $r_i$ that are crescent along the direction $d$ and, therefore, does not exist $\lambda>0$ satisfying~(\ref{eq:desig1}). 
On such cases the index of function $r_i$ belongs to the set $X_0$ and occurs $r_i(\lambda) > p_i(x^*)$. The behavior of function $r_i$ for items 
1.(a), 1.(b) and 2.(a) is illustrated in Figure~\ref{fig:ElementosEmX0eX2eRaiz}, where appears a convex parable and a growing straight line 
for $\lambda>0$.

\begin{figure}[h]
\begin{center}
\begin{tikzpicture}

\draw [->](1,0.5) -- (1,5) node [above]{\small $r_i(\lambda)$};
\draw [->](.5,1) -- (7,1) node [right]{\small $\lambda$};

\draw [dashed](.9,2) -- (7,2) (0.5,2) node {{\tiny$p_{4} (x^*)$}}; \draw [line width=.7pt] (1,2) parabola[bend pos=0.5] bend +(0,4) +(2,-0.8);\draw [->, line width=1.2pt](2.9,2) -- (3.3,2);\draw [dashed](2.9,2) -- (2.9,0.7) node[below] {{\scriptsize $\lambda_{2}^{d}$}};\draw [|->, line width=1.2pt](2.89,1) -- (3.5,1);

\draw [dashed](.9,2.5) -- (7,2.5) (0.5,2.5) node {{\tiny $\textcolor{black}{p_{3} (x^*)}$}};\draw[color=black, line width=.7pt] (1,2.5) sin (3,1.5) cos (6,5);\draw [color=black,<-, line width=1.2pt](4.05,2.5) -- (4.45,2.5);\draw [dashed](4.45,2.5) -- (4.45,.7) node[below] {{\scriptsize \textcolor{black}{$\lambda_1^{d}$}}};\draw [color=black,<-|, line width=1.2pt](3.85,1) -- (4.46,1);

\draw [dashed](.9,3) -- (7,3) (0.5,3) node {{\tiny $\textcolor{black}{p_{2} (x^*)}$}};\draw[color=black, line width=.7pt] (1,3) sin (3.5,1.3) cos (7,5.5);\draw [color=black,<-, line width=1.2pt](5.18,3) -- (5.58,3);

\draw [dashed](.9,3.5) -- (7,3.5) (0.5,3.5) node {{\tiny $p_{1} (x^*)$}};\draw [line width=.7pt] (1,3.5) parabola[bend pos=0.5] bend +(0,1.7) +(2.1,-2);\draw [->, line width=1.2pt](2.6,3.5) -- (3,3.5);

\draw (6.7,3.8) node {$d\in L$};\draw (6.3,0.5) node [color=black] {{\small $\Lambda^{d}\neq\emptyset $}};
\draw [color=black,<->, line width=1.2pt] (2.89,.85) -- (4.46,.85);

\draw [<-, color=black](5.9,4.95) -- (5.8,5.3);
\draw [->, color=black](5.8,5.3) -- (6.68,5);
\draw [color=black] (6,5.5) node {\small $2,3\in X_1$};

\draw [<-](2.45,4.88) -- (2.8,5);
\draw [->](2.8,5) -- (2.3,4.1);
\draw [] (3.1,5.2) node {\small $1,4\in X_2$};

\end{tikzpicture}
\caption{A search direction by $x'$ that dominates $x^*\in$ \textit{Leff}$(VQFP')_{x^*}$, $X_0=\emptyset$, $X_1=\{2,3\}$, $X_2=\{1,4\}$ and $\Lambda^{d}\neq \emptyset$}\label{fig:Direcao2BuscaSemRestricoes}
\end{center}
\end{figure}
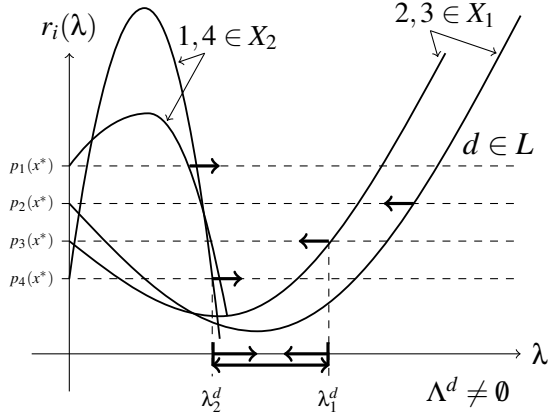

In Figure~\ref{fig:ElementosEmX1eRaiz} is represented the event of item 1.(c) above, a function $r_i$ that is the graph of a 
convex parable. In item 1.(c), exists $\lambda>0$ satisfying~(\ref{eq:desig1}) and it gives the conditions for $r_i$ to be a convex quadratic 
function, in which the equation $r_i(\lambda)-p_i(x^*)=0$ has a positive root $\bar{\lambda}=\frac{-2\nabla p_i(x^*)^T d}{d^T \nabla^2 p_i(x^*)d}$. 
Then, to occur $p_i(x^* +\lambda d) \leqq p_i(x^*)$, we must have $\lambda\in\left(0,\frac{-2\nabla p_i(x^*)^T d}{d^T \nabla^2 p_i(x^*) d}\right]$. 
We exclude zero in this interval because we are interested on points $x^* + \lambda d$ distinct to $x^*$.

Items 2.(b), 2.(c), 3.(b) and 3.(c) represent functions $r_i$  which are not crescent along the unitary direction $d$, therefore, 
all $\lambda\in (0,\infty)$ satisfy~(\ref{eq:desig1}). The behavior of functions $r_i$ for such items are illustrated in 
Figure~\ref{fig:ElementosEmX1eRaiz}. Item 2.(b) imply that for one $i\in I$, the function $r_i$ is constant, because your linear and quadratic term are zero. 
The items 2.(c), 3.(b) and 3.(c) occur when $d^T \nabla^2 p_i(x^*) d < 0$ and $\nabla p_i(x^*)^T d \leqq 0$, or when 
$d^T \nabla^2 p_i(x^*) d \leqq 0$ and $\nabla p_i(x^*)^T d < 0$, for one $i\in I$. In this case, it is possible to occur two situations whenever 
we have $X_0=\emptyset$ and $X_2 \neq \emptyset$ on the direction $d$. In the first, if $X_1=\emptyset$ then $\lambda_1^d=\infty$ and the inequality 
$r_i(\lambda)\leqq p_i(x^*)$ is valid $\lambda>0$ determined by the positive roots of equations $r_i(\lambda)-p_i(x^*)=0$ whose indices
$i \in X_2$. In the second, when $X_1 \neq \emptyset$, if exists one index $k \in I$ such that $d^T \nabla^2 p_k(x^*) d \leqq 0$ and 
$\nabla p_k(x^*)^T d \leqq 0$, the value of $\lambda_1^d$ is not influenced by the function $r_k$, because $\lambda_1^d$ is determined by the 
positive roots of equations $r_i(\lambda)-p_i(x^*)=0$ whose indices $i \in X_1$.

\begin{figure}[h]
\begin{center}
\begin{tikzpicture}

\draw [->](1,0.5) -- (1,5) node [above]{\small $r_i(\lambda)$};
\draw [->](.5,1) -- (7,1) node [right]{\small $\lambda$};

\draw [dashed](.9,2) -- (7,2) (0.5,2) node {{\tiny $p_{4} (x^*)$}}; \draw [line width=.7pt] (1,2) parabola[bend pos=0.5] bend +(0,2) +(2,-0.5);

\draw [dashed](.9,2.5) -- (7,2.5) (0.5,2.5) node {{\tiny $ \textcolor {black}{p_{3} (x^*)}$}};\draw[color=black, line width=.7pt] (1,2.5) sin (2.2,1.5) cos (4.5,5);

\draw [dashed](.9,3) -- (7,3) (0.5,3) node {{\tiny $p_{2} (x^*)$}};\draw [line width=.7pt] (1,3) parabola[bend pos=0.5] bend +(0,2) +(5,-1.5);

\draw [dashed](.9,3.5) -- (7,3.5) (0.5,3.5) node {{\tiny $\textcolor {black}{p_{1} (x^*)}$}};\draw[ line width=1.2pt, color=black] (1,3.5) cos (4,5);

\draw (6.3,4) node {$d \notin L$};\draw (6.3,0.6) node {{\small $X_0 \neq \emptyset $}};

\draw [->, color=black](3,5) -- (3.57,4.77);
\draw [color=black] (2.5,5) node {\small $1\in X_0$};

\draw [<-, color=black](4.5,4.8) -- (5.2,5);
\draw [color=black] (5.5,5.2) node {\small $3\in X_1$};

\draw [<-](3,1.7) -- (4.3,1.4);
\draw [->](4.3,1.4) -- (5.8,1.8);
\draw [] (4.5,1.2) node {\small $2,4\in X_2$};

\end{tikzpicture}
\caption{A search direction by $x'$ that dominates $x^*\in$ \textit{Leff}$(VQFP')_{x^*}$, $X_0=\{1\}$, $X_1=\{3\}$ and $X_2=\{2,4\}$}\label{fig:Direcao3BuscaSemRestricoes}
\end{center}
\end{figure}
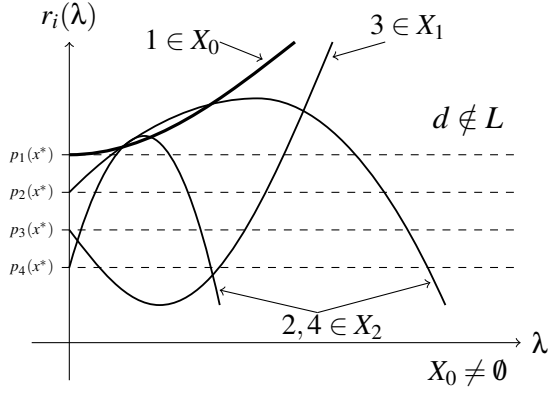

In Figure~\ref{fig:ElementosEmX0eX2eRaiz} is represented what occurs in item 3.(a). There, the function $r_i$ is the graphic 
of a concave parabola. In this item, exists $\lambda>0$ satisfying~(\ref{eq:desig1}) and equation $r_i(\lambda)-p_i(x^*)=0$ has a positive root 
$\bar{\lambda}=\frac{-2\nabla p_i(x^*)^T d}{d^T \nabla^2p_i(x^*)d}$. Therefore, to occur $p_i(x^* +\lambda d) \leqq p_i(x^*)$ we must have 
$\lambda\in\left[\frac{-2\nabla p_i(x^*)^T d}{d^T \nabla^2 p_i(x^*) d},\infty\right)$.

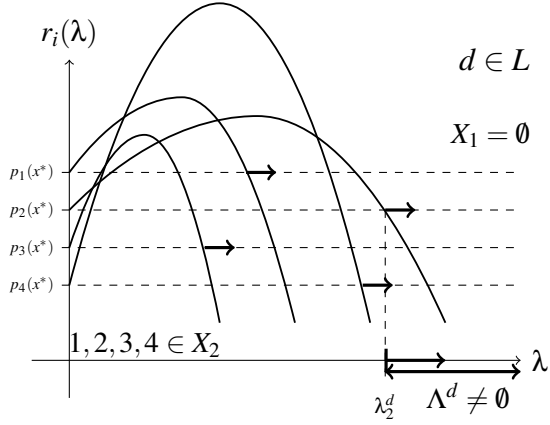
\begin{figure}[h]
\begin{center}
\begin{tikzpicture}

\draw [->](1,0.5) -- (1,5) node [above]{\small $r_i(\lambda)$};
\draw [->](.5,1) -- (7,1) node [right]{\small $\lambda$};

\draw [dashed](.9,2) -- (7,2) (0.5,2) node {{\tiny$p_{4} (x^*)$}}; \draw [line width=.7pt] (1,2) parabola[bend pos=0.5] bend +(0,4) +(4,-0.5);\draw [->, line width=1.2pt](4.9,2) -- (5.3,2);

\draw [dashed](.9,2.5) -- (7,2.5) (0.5,2.5) node {{\tiny$ p_{3} (x^*)$}};\draw [line width=.7pt] (1,2.5) parabola[bend pos=0.5] bend +(0,2) +(2,-1); \draw [->, line width=1.2pt](2.8,2.5) -- (3.2,2.5);

\draw [dashed](.9,3) -- (7,3) (0.5,3) node {{\tiny$p_{2} (x^*)$}};\draw [line width=.7pt] (1,3) parabola[bend pos=0.5] bend +(0,2) +(5,-1.5);\draw [->, line width=1.2pt](5.2,3) -- (5.6,3);\draw [dashed](5.2,3) -- (5.2,0.7) node[below] {{\scriptsize $\lambda_{2}^{d}$}};\draw [|->, line width=1.2pt](5.19,1) -- (6,1);

\draw [dashed](.9,3.5) -- (7,3.5) (0.5,3.5) node {{\tiny$p_{1} (x^*)$}};\draw [line width=.7pt] (1,3.5) parabola[bend pos=0.5] bend +(0,2) +(3,-2);\draw [->, line width=1.2pt](3.35,3.5) -- (3.75,3.5);

\draw (6.65,5) node {$d\in L$};\draw (6.3,0.5) node [color=black] {{\small $\Lambda^{d}\neq\emptyset$}};\draw (6.6,4) node {{\small \textcolor{black}{$X_1=\emptyset $}}};
\draw [color=black,<->, line width=1.2pt] (5.19,.85) -- (7,.85);

\draw [] (2,1.2) node {\small$1,2,3,4\in X_2$};

\end{tikzpicture}

\caption{A search direction by $x'$ that dominates $x^*\in$ \textit{Leff}$(VQFP')_{x^*}$, $X_0=\emptyset$, $X_1=\emptyset$, $X_2=\{1,2,3,4\}$ and $\Lambda^{d}\neq\emptyset$}\label{fig:Direcao4BuscaSemRestricoes}
\end{center}
\end{figure}

Given an unitary direction $d$, we look for a nonempty interval $\Lambda^d$ such that $x'=x^* +\lambda d$, 
$p_i(x')\leqq p_i(x^*)$, $\forall i\in I$, where $\lambda\in \Lambda^d$ and $X_0=\emptyset$. Hence, we must verify the items above for $m$ 
functions. In Fugures~\ref{fig:Direcao1BuscaSemRestricoes},~\ref{fig:Direcao2BuscaSemRestricoes},~\ref{fig:Direcao3BuscaSemRestricoes}
and~\ref{fig:Direcao4BuscaSemRestricoes} are plotted the points of the parabolas that are the graphics of functions 
$r_i(\lambda)=p_i(x^* +\lambda d)$, $i\in\{1,2,3,4\}$, at coordinates $(\lambda,r_i(\lambda))$. Figure~\ref{fig:Direcao2BuscaSemRestricoes} illustrates 
a case which $\Lambda^d\neq \emptyset$, $X_1=\{2,3\}$ and $X_2=\{1,4\}$, where we see that $\lambda_2^d \leqq \lambda_1^d$. Note also that 
if $\lambda\in \Lambda^d$ and $x'=x^* +\lambda d$, we obtain $p_i(x')\leqq p_i(x^*)$, $\forall i\in\{1,2,3,4\}$.

Figure~\ref{fig:Direcao3BuscaSemRestricoes} illustrates a case which $X_0 \neq \emptyset$, where we have $X_0=\{1\}$, 
$X_1=\{3\}$, and $X_2=\{2,4\}$. In this case, $r_1(\lambda)$ always increases for $\lambda > 0$, then $x^*$ is globally Pareto-efficient 
solution in this direction, that is, whenever occurs $X_0 \neq \emptyset$ in a direction $d$, the solution $x^*$ is globally Pareto-efficient 
solution in this direction. On the other hand, given a direction $d$, if $X_0 = \emptyset$ and $X_1 \neq \emptyset$, then not occur $X_2=\emptyset$ 
in this direction, because the solution $x^*$ would not be locally Pareto-efficient solution, since $p(x^*+\lambda d)\leq p(x^*)$ would be valid 
for all $\lambda>0$. This situation can be observed in Figures~\ref{fig:Direcao1BuscaSemRestricoes} and~\ref{fig:Direcao2BuscaSemRestricoes}, 
needing only to assume $X_2=\emptyset$ in both figures. In fact, suppose that $X_2=\emptyset$ in figure~\ref{fig:Direcao2BuscaSemRestricoes}, 
then when ignoring the graphics of functions $r_1$ and $r_4$, we obtain $r_1(\lambda)\leqq p_1(x^*)$, $r_2(\lambda)< p_2(x^*)$, 
$r_3(\lambda) < p_3(x^*)$ and $r_4(\lambda)\leqq p_4(x^*)$ for all $\lambda>0$. That is, $p(x^*+\lambda d)\leq p(x^*)$ is valid for all $\lambda>0$.

Therefore, whenever occur $X_0 = \emptyset$ in a determined direction $d$ and $x^*$ being locally Pareto-efficient optimal solution, two 
important situations to our analysis may happen: $X_1\neq \emptyset$ and $X_2\neq \emptyset$, or $X_1= \emptyset$ and $X_2\neq \emptyset$. Figure~\ref{fig:Direcao4BuscaSemRestricoes} 
illustrates the case $X_0=\emptyset$, $X_1=\emptyset$ and only $X_2\neq \emptyset$, where $X_2=\{1,2,3,4\}$, $X_0=\emptyset$ and $X_1=\emptyset$, 
then we obtain $\lambda_1^d=\infty$, $\lambda_2^d<\lambda_1^d$ and $\Lambda^d\neq \emptyset$. Note in the unconstrained problems that, in the 
direction $d$, always is possible calculation the positive root of the equation $r_i(\lambda)-p_i(x^*)=0$, $i\in X_2$, and if occur $X_0=\emptyset$, 
$X_1=\emptyset$ and $X_2\subseteq \{1,\ldots,m\}$ with $|X_2|\geqq 1$ (if not occur item 2.(b), we obtain $X_2=\{1,\ldots,m\}$), then we 
always have $\Lambda^d \neq \emptyset$ in this direction and it is possible to find a point $x'$ that dominates $x^*$.

We can conclude that if $X_0= \emptyset$, $X_1\neq \emptyset$ and $X_2\neq \emptyset$ in a direction $d$, the values of 
$\lambda >0$, such that, $ \frac{f_i(x^*+\lambda d)}{g_i(x^*+\lambda d)}\leqq \frac{f_i(x^*)}{g_i(x^*)}$, $\forall i\in I$, must satisfy

\begin{equation}\label{eq:IdentidadeLambdaInicialMenor}
\lambda \leqq \ \frac{-2\nabla p_i(x^*)^T d}{d^T \nabla^2 p_i(x^*) d}\ \ \mbox{if}\ \ i\in X_1\quad\ \ \mbox{and}\quad\ \ \lambda \geqq \ \frac{-2\nabla p_i(x^*)^T d}{d^T \nabla^2 p_i(x^*) d}\ \ \mbox{if}\ \ i\in X_2.
\end{equation} 
Therefore, to grant that $\frac{f(x^* + \lambda d)}{g(x^* + \lambda d)}\leq \frac{f(x^*)}{g(x^*)}$ have solution 
in the unitary direction $d$, we must choose correctly $\lambda>0$ such that $\lambda_2^d \leqq \lambda \leqq \lambda_1^d$, and this 
choose must grant at least that one inequality be strict. 

Next, some possibility which harden the demonstration of our results are explained. For 
example, we must see what happens when there exists a direction $d$ such that $X_0= \emptyset$, $X_1\neq \emptyset$, $X_2\neq \emptyset$ and 
$\lambda_2^d=\lambda_1^d=\frac{-2\nabla p_i(x^*)^T d}{d^T \nabla^2 p_i(x^*) d}$, for all $i\in X_1 \cup X_2$, and also when there exists a 
direction $d$ such that $X_0= \emptyset$, $X_1= \emptyset$ and $X_2= \emptyset$. First, suppose that we obtain $X_0= \emptyset$, $X_1\neq \emptyset$ and $X_2\neq \emptyset$, then:

\begin{enumerate}
\item[4.] if $|X_1|\geqq 2$,
\begin{enumerate}
\item $|X_2|\geqq 1$ and $ \lambda_2^d=\lambda_1^d=\frac{-2\nabla p_i(x^*)^T d}{d^T \nabla^2 p_i(x^*) d}$, for all $i\in X_1 \cup X_2$. 
As $\Lambda_1^d=(0,\lambda_1^d]$, we obtain $\Lambda^d=\Lambda_2^d\cap\Lambda_1^d=\{\lambda_1^d\}$, and then exists $\lambda\in\Lambda^d$. However, 
it is elementary to verify that $p(x^*+\lambda d)\leq p(x^*)$ does not have solution for each $\lambda>0$. 
We need to be careful in this case, because we obtain a nonempty interval $\Lambda^d=\{\lambda_1^d\}$, however $p(x^*+\lambda_1^d d)\leq p(x^*)$ does not have solution.

\item $|X_2|\geqq 1$ and there exists $j\in X_2$ such that $\frac{-2\nabla p_j(x^*)^T d}{d^T \nabla^2 p_j(x^*) d}< \lambda_2^d=\lambda_1^d$. 
As $\Lambda_1^d=(0,\lambda_1^d]$, we obtain $\Lambda^d=\Lambda_2^d\cap\Lambda_1^d=\{\lambda_1^d\}$. Hence, there exists 
$\lambda\in\Lambda^d$ and as $j\in X_2$ such that $ \frac{-2\nabla p_j(x^*)^T d}{d^T \nabla^2 p_j(x^*)d} < \lambda_2^d $, it 
is elementary to verify that, unlike item 4.(a), $p(x^*+\lambda_2^d d)\leq p(x^*)$ has solution. 

\item $|X_2|\geqq 1$ and, analogously item 4.(b), there exists $j\in X_1$ such that $ \lambda_2^d=\lambda_1^d < \frac{-2\nabla p_j(x^*)^T d}{d^T \nabla^2 p_j(x^*) d}$. Hence, 
we obtain $\Lambda^d=\Lambda_2^d \cap \Lambda_1^d=\{\lambda_1^d\}$, there exists $\lambda\in\Lambda^d$ and as $j\in X_1$ such that 
$ \lambda_1^d < \frac{-2\nabla p_j(x^*)^T d}{d^T \nabla^2 p_j(x^*)d}$, it is elementary to verify that, unlike item 4.(a), $p(x^*+\lambda_1^d d)\leq p(x^*)$ has 
solution. 

\end{enumerate}
\end{enumerate}

Figure~\ref{fig:Direcao2BuscaSemRestricoes} helps to understand the items discussed above. Now, suppose that given a direction $d$ we obtain 
$X_0=\emptyset$, $X_1=\emptyset$ and $X_2=\emptyset$. This means that for the existence of positive solutions to inequality~(\ref{eq:desig1}), 
only item 2.(b) occurs for all $i\in I$. That is, simultaneosly $\nabla p_i(x^*)^T d=d^T\nabla^2 p_i(x^*) d=0$, for all $i\in I$, and inequality~(\ref{eq:desig1}) is valid for 
all $\lambda\in (0,\infty)$. In this case, among other hypothesis, the following can occurs.

\begin{enumerate}
\item[5.] $d$ is orthogonal to set $\left\{\nabla p_i(x^*)\right\}_{i\in I}$, 
\begin{enumerate}
\item and $d$ is orthogonal to set $\left\{\nabla^2 p_i(x^*) d \right\}_{i\in I}$,

\item and there exists $\bar{I}\subseteq I$, such that, $d$ is orthogonal to set $\left\{\nabla^2 p_i(x^*) d\right\}_{i\in \bar{I}}$ and 
$d\in Nu(\nabla^2 p_i(x^*))$, for all $i\in I \setminus \bar{I}$, where $Nu(\nabla^2 p_i(x^*))$ is the Kernel of the matrix $\nabla^2 p_i(x^*)$.
\end{enumerate}

\item[6.] There exists $\bar{I}\subseteq I$, such that, $d$ is orthogonal to set $\left\{\nabla p_i(x^*)\right\}_{i\in \bar{I}}$ and $\nabla p_i(x^*)=0$, for all $i\in I \setminus \bar{I}$, 
\begin{enumerate}
\item and $d$ is orthogonal to set $\left\{\nabla^2 p_i(x^*) d \right\}_{i\in I}$,

\item and there exists $\tilde{I}\subseteq I$, such that, $d$ is orthogonal to set $\left\{\nabla^2 p_i(x^*) d\right\}_{i\in \tilde{I}}$ and 
$d\in Nu(\nabla^2 p_i(x^*))$, for all $i\in I \setminus \tilde{I}$.
\end{enumerate}
\end{enumerate}

We believe that such events 4.(a)--(c), 5.(a)--(b) or 6.(a)--(b) are rare, because in general the objective functions of the vector 
optimization problems are conflicting. Aiming to facilitate the demonstrations of our results, we impose the following two conditions.

{\tcon{Let $d\in \partial B(0,1)$ and suppose that $X_0=\emptyset$, $X_1\neq\emptyset$ and $X_2\neq\emptyset$. If 
$ \lambda_2^d=\lambda_1^d=\frac{-2\nabla p_i(x^*)^T d}{d^T \nabla^2 p_i(x^*) d}$, for all $i\in X_1 \cup X_2$, then we define $\Lambda^d:=\emptyset$.}\label{cond:LambdaTudoIgual}}

{\tcon{Let $x^*\in\mbox{\textit{Leff}}(VQFP')_{x^*}$ and $d\in \partial B(0,1)$. Then does not simultaneosly occurs 
$\nabla p_i(x^*)^T d=d^T\nabla^2 p_i(x^*) d=0$, for all $i\in I$.}\label{cond:TudoNulo}}\\

Note that if we were dealing with weakly Pareto-efficient optimal solutions, items 4--6 and Conditions~\ref{cond:LambdaTudoIgual} and~\ref{cond:TudoNulo} would not be necessary.

{\tlem{Let $x^*\in\mbox{\textit{Leff}}(VQFP')$. Then $x'$ dominates $x^*$ if and only if there exists $d\in \partial B(0,1)$ such that 
$X_0=\emptyset$, $X_2\neq \emptyset$ and $\Lambda^d\neq \emptyset$. In this case, there exists $\lambda^* \in \Lambda^d$ such that $x'=x^* + \lambda^* d$.}\label{lem:LambdaVazioNovaVersao}}\\

{\prova{$(\Rightarrow)$ If $x^*\in\mbox{\textit{Leff}}(VQFP')$, then from Theorem~\ref{teo:ProblemaAssociado} 
$x^*\in\mbox{\textit{Leff}}(VQFP')_{x^*}$ and there exists a neighborhood $N(x^*)$, such that, $ f(x)-\frac{f(x^*)}{g(x^*)}g(x)\leq 0$ does not have 
solution for each $x\in N(x^*)$, and by~(\ref{eq:DefineFuncoesP}), $p(x)\leq p(x^*)$ does not have solution for each $x\in N(x^*)$. Consider 
$N(x^*)$ as a ball of $\R^n$ whith radius $\bar{\lambda}$, that is, $p(x)\leq p(x^*)$ does not have solution for each 
$x\in N(x^*)\equiv B(x^*,\bar{\lambda})$ and, if $d\in \partial B(0,1)$, $p(x)\leq p(x^*)$ does not have solution for each $x=x^*+\lambda d$ 
such that $\|x-x^*\|\leqq \lambda$, for all $\lambda\in (0,\bar{\lambda})$, where $\|\cdot\|$ is the Euclidean norm. Therefore, if $x'$ dominates 
$x^*$, $x'\notin B(x^*,\bar{\lambda})$ and there exists $d\in \partial B(0,1)$ such that $x'=x^* + \lambda^* d$ and $\lambda^* \geqq \bar{\lambda}$. 
We need to verify what happens to the sets $X_0$, $X_1$, $X_2$ and $\Lambda^d$ in this direction $d$. Figure~\ref{fig:VizinhancaLema31} 
illustrates an unidimensional neighborhood $B(0,\bar{\lambda})$ at point $\lambda=0$, which represents an $n$-dimensional neighborhood 
$B(x^*,\bar{\lambda})$ at point $x^*$. We observe that in a fixed direction $d$, used as example, there exists functions $r_i$ whose graphics 
show that their indices $i$ belong to $X_0$ (increasing curve in dashed line), or $i$ belong to $X_1$ (convex parabola in continous line), 
or $i$ belong to $X_2$ (concave parabola in dashed line), and more, whenever $\lambda$ satisfy items 3.(b) and 3.(c) for any index $i$ 
(decreasing curve in continous line). It is also possible to observe that if $X_0=\emptyset$ in this figure, $r_i(\lambda)\leqq p_i(x^*)$, for all $i\in I$, 
only have solution if $\lambda \in \Lambda^d$ (represented below of the $\lambda$-axis of the cartesian plane). That is, if $X_0=\emptyset$ in 
figure~\ref{fig:VizinhancaLema31}, $p(x^*+\lambda d)\leq p(x^*)$ has solution if and only if $\lambda_2^d \leqq\lambda\leqq \lambda_1^d$. 
The set $\Lambda_1^d$ is represented above the $\lambda$-axis, while the set $\Lambda_2^d$ is represented above the straight line $p_i(x^*)$ in this figure.

Thus, given an arbitrary $d\in \partial B(0,1)$, suppose that $X_0\neq \emptyset$ in this direction. Then there exists an index $i\in I$, 
such that, the item 1.(a), 1.(b) or 2.(a) is satisfied, that is, $d^T \nabla^2 p_i(x^*) d>0$ and $\nabla p_i(x^*)^T d \geqq 0$ or $d^T \nabla^2 p_i(x^*) d\geqq 0$ 
and $\nabla p_i(x^*)^T d>0$. In this case, $r_i(\lambda)$ grows indenifinitely for $\lambda>0$ and therefore given any heighborhood of $x^*$ in 
direction $d$ (e.g., $B(0,\lambda)$ with $\lambda\geqq\bar{\lambda}>0$, in Figure~\ref{fig:VizinhancaLema31}), $p(x)\leq p(x^*)$ does not have 
solution for each $x=x^*+\lambda d$ in this neighborhood, hence does not exist $x'=x^* + \lambda d$ such that $x'$ dominates $x^*$ in this 
direction. Therefore, $X_0$ has to be empty.

\begin{figure}[h]
\begin{center}
\begin{tikzpicture}


\draw [latex-), black!70, line width=1.8pt] (1,3) -- (2.1,3);
\draw [white!95, line width=5pt] (1,3.1) -- (1.5,3.1);
\draw [white!95, line width=5pt] (1,2.9) -- (1.5,2.9);

\draw [->](1,0.5) -- (1,5) node [above]{\small $r_i(\lambda)$};
\draw [->](.5,1) -- (7,1) node [right]{\small $\lambda$};

\draw[color=black, line width=.7pt] (1,3) sin (3,1.5) cos (6.6,5);
\draw [dashed](5.2,3) -- (5.2,.7);
\draw [line width=1pt, color=white](5.2,2) -- (5.2,2.2);
\draw (5.03,1) node[above] {{\scriptsize \textcolor{black}{$\lambda_1^{d}$}}};
\draw [color=black,<-|, line width=1.2pt] (1,1) -- (5.2,1);
\draw (1.4,1) node[above] {{\scriptsize \textcolor{black}{$\Lambda_1^{d}$}}};
\draw [color=black,<->, line width=1.2pt] (3.6,.85) -- (5.2,.85);\draw (4.4,.85) node[below] {{\scriptsize \textcolor{black}{$\Lambda^{d}$}}};

\draw [dashed](1,3) -- (7,3) (0.5,3) node {{\scriptsize $p_{i} (x^*)$}};
\draw [dashed, line width=.7pt] (1,3) parabola[bend pos=0.5] bend +(0,2.5) +(3,-1.4);
\draw [|->, line width=1.2pt](3.58,3) -- (7,3);
\draw (6.95,3.05) node[right] {{\scriptsize \textcolor{black}{$\Lambda_2^{d}$}}};
\draw [dashed](3.6,3) -- (3.6,0.7);
\draw (3.78,3) node[above] {{\scriptsize $\lambda_{2}^{d}$}};


\draw [<-](2.94,4.6) -- (3.4,5);
\draw (3.5,5.2) node {\small $i\in X_2$};

\draw [color=black,->](4.9,5) -- (4.1,4.8);
\draw (5,5.2) node [color=black]{\small $i\in X_0$};

\draw [color=black,->](6.5,4) -- (6.2,4.3);
\draw (6.57,3.8) node [color=black]{\small $i\in X_1$};

\draw[dashed, color=black, line width=.7pt] (1,3) cos (4.2,5);
\draw[color=black, line width=.7pt] (1,3) cos (4.15,1.2);

\draw (.5,2.05) node {{\scriptsize $B(0,\bar{\lambda})$}};
\draw [->] (.55,2.20) -- (1.9,2.93);
\draw [->, black!80] (1,3) -- (2.08,3);
\draw (2.2,3.15) node {{\scriptsize $\bar{\lambda}$}};

\draw [color=black,->](4.7,2) -- (4.05,1.35);
\draw [color=black] (5.5,2.1) node {\tiny items 3(b) and 3(c)};

\draw [black!70, line width=1.3pt] (1,3) -- (2.1,3);

\end{tikzpicture}
\caption{Neighborhood of point $x^*$ and the behavior of the functions $r_i$}\label{fig:VizinhancaLema31}
\end{center}
\end{figure}

Suppose now that $X_0= \emptyset$ and $X_2=\emptyset$ in direction $d$. As $x^*\in$ \textit{Leff}$(VQFP')_{x^*}$ in 
$B(x^*,\bar{\lambda})$, $x^*$ is also locally Pareto-efficient optimal solution in $B(0,\bar{\lambda})$ in direction $d$. Then in this 
direction cannot occur $X_1\neq\emptyset$ and none of items 2.(c), 3.(b) or 3.(c) (decreasing curve in continous line shown in 
Figure~\ref{fig:VizinhancaLema31}) can be satisfied, that is, cannot occur $d^T \nabla^2 p_i(x^*) d > 0$ and $\nabla p_i(x^*)^T d <0$, 
as well, cannot occur $d^T \nabla^2 p_i(x^*) d \leqq 0$ and $\nabla p_i(x^*)^T d < 0$, or $d^T \nabla^2 p_i(x^*) d < 0$ and 
$\nabla p_i(x^*)^T d \leqq 0$. Hence, only item 2.(b) could be satisfied in this direction, that is, $d^T \nabla^2 p_i(x^*) d = 0$ 
and $\nabla p_i(x^*)^T d = 0$, $\forall i\in I$. Wich is impossible, due to Condition~\ref{cond:TudoNulo}. Therefore, $X_0= \emptyset$ and $X_2\neq\emptyset$ in direction $d$.

Suppose now that $X_0= \emptyset$ and $X_2\neq\emptyset$, but $\Lambda^d=\emptyset$ in direction $d$. With such hypothesis, 
$X_1=\emptyset$ cannot occurs in direction $d$, otherwise we would have $\lambda_1^d=\infty$ and $\Lambda^d\neq\emptyset$. Then, 
if $X_0= \emptyset$, $X_1\neq \emptyset$, $X_2\neq\emptyset$ and $\Lambda^d=\emptyset$ two cases may occur. Either $\lambda_1^d<\lambda_2^d$ 
and $p(x^*+\lambda d)\leq p(x^*)$ does not have solution for $\lambda>0$, or by Condition~\ref{cond:LambdaTudoIgual} and item 4.(a), 
$ \lambda_2^d=\lambda_1^d=\frac{-2\nabla p_i(x^*)^T d}{d^T \nabla^2 p_i(x^*) d}$, for all $i\in X_1 \cup X_2$, and $p(x^*+\lambda d)\leq p(x^*)$ 
does not have solution for $\lambda>0$. Therefore, $x'=x^*+\lambda d$ does not dominates $x^*$ in this direction. We conclude that if $x'$ dominates 
$x^*$, then there exists $d\in \partial B(0,1)$ where $X_0= \emptyset$ and $X_2\neq\emptyset$, and more $\lambda^*\in\Lambda^d\neq\emptyset$ 
such that $x'=x^*+\lambda^* d$.

$(\Leftarrow)$ Suppose that in the direction $d$, we obtain $X_0= \emptyset$, $X_2\neq\emptyset$ and $\Lambda^d\neq\emptyset$. If $X_1=\emptyset$, 
then there exists $\lambda^*\in \Lambda^d$ (see Figure~\ref{fig:Direcao4BuscaSemRestricoes}), such that, if $x'=x^*+\lambda^*d$ then $p(x')\leq p(x^*)$ has solution. 
If $X_1\neq\emptyset$, then either items 4.(b)--(c) grant that exists $\lambda^*\in \Lambda^d$, such that, if $x'=x^*+\lambda^*d$ then 
$p(x')\leq p(x^*)$ has solution, or $\lambda_2^d < \lambda_1^d$ and there exists $\lambda^*\in \Lambda^d=[\lambda_2^d, \lambda_1^d]$, such that, if $x'=x^*+\lambda^*d$ 
then $p(x')\leq p(x^*)$ has solution. Therefore, $x'$ dominates $x^*$. }}

{\tteo{Let $x^*\in$ \textit{Leff}$(VQFP')$. Then $x^*\in$ \textit{Eff}$(VQFP')$ if and only if for all $d\in \partial B(0,1)$, $X_0\neq \emptyset$ 
or $\Lambda^d =\emptyset$ whenever $X_1\neq \emptyset$.}\label{teo:LambdaVazio}}\\

{\prova{$(\Rightarrow)$ Suppose that exists $d\in \partial B(0,1)$ where $X_1\neq \emptyset$, $X_0=\emptyset$ and $\Lambda^d\neq \emptyset$. 
As $x^*\in$ \textit{Leff}$(VQFP')$, from Theorem~\ref{teo:ProblemaAssociado} $x^*\in$ \textit{Leff}$(VQFP')_{x^*}$ 
and $X_2\neq \emptyset$ in this direction. By Lemma~\ref{lem:LambdaVazioNovaVersao}, there exists $\lambda^*\in \Lambda^d$ and 
$x'=x^* +\lambda^* d$ such that $p(x')\leq p(x^*)$ has solution. By~(\ref{eq:DefineFuncoesP}), 
$ \frac{f(x')}{g(x')}\leq \frac{f(x^*)}{g(x^*)}$ has solution. What contradicts $x^*\in$ \textit{Eff}$(VQFP').$

$(\Leftarrow)$ By Theorem~\ref{teo:ProblemaAssociado}, $x^*\in$ \textit{Leff}$(VQFP')_{x^*}$. Let an arbitrary $d\in \partial B(0,1)$ and suppose 
that $X_0 \neq \emptyset$ in direction $d$. Then there exists an index $i\in I$, such that, one of items 1.(a) and 1.(b) or 2.(a) is satisfied, 
that is, $d^T \nabla^2 p_i(x^*) d>0$ and $\nabla p_i(x^*)^T d \geqq 0$, or $d^T \nabla^2 p_i(x^*) d\geqq 0$ and $\nabla p_i(x^*)^T d > 0$. Hence, 
$r_i(\lambda)$ grows indefinitely for $\lambda>0$ and given any neighborhood of $x^*$ in direction $d$,  as large as it is, 
$p(x')\leq p(x^*)$ does not have solution for each $x'=x^*+\lambda d$ in this neighborhood. Therefore, $x'=x^* + \lambda d$ does not dominates 
$x^*$ and $x^*$ is Pareto-efficient optimal solution in this direction. Figure~\ref{fig:Direcao3BuscaSemRestricoes} illustrates this possibility, 
where the function $r_1$, $1\in X_0$, grows indefinitely for $\lambda>0$. On the other hand, suppose that $X_1\neq \emptyset$ and $X_0 = \emptyset$, but $\Lambda^d=\emptyset$ in direction $d$. 
As $x^*\in$ \textit{Leff}$(VQFP')_{x^*}$, necessarily $X_2\neq \emptyset$ in this direction. We have then $X_0= \emptyset$, 
$X_1\neq \emptyset$, $X_2\neq\emptyset$, $\Lambda^d=\emptyset$ and two cases may occur. Either $\lambda_1^d<\lambda_2^d$ and 
$p(x^*+\lambda d)\leq p(x^*)$ does not have solution for each $\lambda>0$, or by Condition~\ref{cond:LambdaTudoIgual} and item 4.(a), 
$\lambda_2^d=\lambda_1^d=\frac{-2\nabla p_i(x^*)^T d}{d^T \nabla^2 p_i(x^*) d}$, for all $i\in X_1\cup X_2$, and $p(x^*+\lambda d)\leq p(x^*)$ 
does not have solution for each $\lambda>0$. Therefore, $x'=x^* + \lambda d$ does not dominates $x^*$ and $x^*$ is Pareto-efficient optimal solution in this direction. Figure~\ref{fig:Direcao1BuscaSemRestricoes} 
illustrates this possibility, there we observe $X_0= \emptyset$, $X_1=\{1,3\}$, $X_2=\{2,4\}$, $\lambda_1^d < \lambda_2^d$, 
$\Lambda^d=\emptyset$ and $p(x^* + \lambda d)\leq p(x^*)$ does not have solution for each $\lambda>0$. 

Therefore, given an arbitrary direction $d\in \partial B_n(0,1)$, if $X_0\neq \emptyset$ we obtain $x^*$ non-dominated, or 
whenever $X_1\neq \emptyset$ and $\Lambda^d =\emptyset$ we obtain again $x^*$ non-dominated. Hence, does not exist another point 
$x'=x^*+\lambda d$ that dominates $x^*$, and $x^* \in$ \textit{Eff}$(VQFP')$. }}\\

Given a direction $d$, in the proof of Theorem~\ref{teo:LambdaVazio} when $X_0=\emptyset$ and $\Lambda^d=\emptyset$, we require 
that $X_1$ be nonempty, otherwise the Condition~\ref{cond:TudoNulo} grant us only $X_2\neq \emptyset$ 
(when $X_0 = \emptyset$ and $X_1 = \emptyset$) and we obtain $\lambda_1^d=\infty$ and $\lambda_2^d < \lambda_1^d$. When this last case occurs in the 
unconstrained problem, we always have $\Lambda^d \neq\emptyset$, that is, there exists $\lambda \in \Lambda^d$ such that 
$x^*+\lambda d$ dominates $x^*$ in direction $d$. Figure~\ref{fig:Direcao4BuscaSemRestricoes} illustrates this possibility, there we observe 
that $X_0=\emptyset$, $X_1=\emptyset$, $X_2=\{1,2,3,4\}$ and the interval $\Lambda^d$ is nonempty. But, this does not always happens 
in the constrained problem, because the feasible set can be limited in this direction.

Note, unlike the results obtained by Osuna-G\'omez et al.~\cite{Art:Osuna}, that on the objective functions  of the auxiliary 
problem (VQFP$'$)$_{x^*}$ we do not require any kind of generalized convexity to obtain the Pareto optimality conditions of Theorem~\ref{teo:LambdaVazio}.

In the results that follows we define $L=\{d\in \partial B(0,1)\ |\ \Lambda^d \neq \emptyset\}$.

{\tcor{Let $x^*\in$ \textit{Leff}$(VQFP')$ and $\beta = \inf\limits_{d\in L}\ \left\{\lambda_{2}^d\right\}$. Then does not exist another 
point $x'\in B(x^*,\beta)$ such that $\frac{f(x')}{g(x')} \leq \frac{f(x^*)}{g(x^*)}$.}\label{cor:LambdaNaoVazio}}\\

{\prova{It follows immediately from Theorem~\ref{teo:LambdaVazio}. If $d\in L$, then there exists $\lambda>0$ in $\Lambda^d$ such that 
another point $x'\in \R^n$ dominates $x^*$ in this direction. However, the first point $x'$ that dominates $x^*$ along the direction $d$ 
is $x'=x^*+\lambda_{2}^d d$. By checking all set $L$, we conclude that the point $x'$ which dominates $x^*$ is $x'=x^*+\beta d$. 
Therefore $x^*$ is $\beta$-efficient and does not exist another point $x'\in B(x^*,\beta)$ such that 
$\frac{f(x')}{g(x')} \leq \frac{f(x^*)}{g(x^*)}$. Figure~\ref{fig:Direcao2BuscaSemRestricoes} illustrates this possibility, there 
we observe that $X_0= \emptyset$, $X_1=\{2,3\}$, $X_2=\{1,4\}$, $\Lambda^d \neq \emptyset$ and if 
$x'=x^*+\lambda_{2}^d d$, $p(x')\leq p(x^*)$ has solution. }}\\

Although complicate your proof, it is possible to rewrite Corollary~\ref{cor:LambdaNaoVazio} to replace the set 
$L$ with the set $L'=\{d \in \partial B_n(0,1)\ |\ X_0=\emptyset\}$ (see in the charts in Figures~\ref{fig:Direcao1BuscaSemRestricoes} 
and~\ref{fig:Direcao2BuscaSemRestricoes}). Next, a very useful result which gives a lower bound for the radius of efficiency of the 
solution $x^*\in$ \textit{Leff}$(VQFP')$ is presented.

{\tcor{Let $x^*\in$ \textit{Leff}$(VQFP')$ and $\displaystyle F(d) = \max\limits_{i\in X_2}\ \left\{ 2\nabla p_i(x^*)^T d \right\}.$ Suppose that 
exists $\rho\in \mathbf{R}$, such that for all $d\in \partial B(0,1)$ we have  $F(d)\geqq \rho$. Then does not exist another point 
$x'\in B(x^*,\frac{\rho}{-\gamma})$ such that $\frac{f(x')}{g(x')} \leq \frac{f(x^*)}{g(x^*)}$, where $\gamma < 0$, 
$\displaystyle \gamma = \min\limits_{i\in I} \left\{\gamma_i\right\}$ and $\gamma_i$ is the smallest negative eigenvalue of the matrix $\nabla^2 p_i(x^*)$, 
$i\in I$.}\label{cor:FMaiorRho}}\\ 

{\prova{By the proof of theorem~\ref{teo:LambdaVazio} and the value of $\beta$ in Corollary~\ref{cor:LambdaNaoVazio}, if we search a point 
$x'=x^* +\lambda d$ that dominates $x^*$, we have to find a value $\lambda>0$ that satisfy
\begin{eqnarray*}
\lambda \geqq \beta = \inf\limits_{d\in L} \ \left\{\max\limits_{i \in X_2}\ \left\{\ \frac{-2\nabla p_i(x^*)^T d}{d^T \nabla^2 p_i(x^*) d}\ \right\}\right\}.
\end{eqnarray*}
Given the symmetric and diagonalizable matrix $\nabla^2 p_i(x^*)$, $i\in X_2$, we obtain $0 > d^T \nabla^2 p_i(x^*) d \geq \gamma_i$, 
where $\gamma_i$ is the smallest negative eigenvalue of the matrix $\nabla^2 p_i(x^*)$, therefore

$$\begin{array}{rcl}
\lambda\ \geqq\ \beta &=& \inf\limits_{d\in L}\ \left\{\max\limits_{i \in X_2}\ \left\{\ \displaystyle \frac{-2\nabla p_i(x^*)^T d}{d^T \nabla^2 p_i(x^*) d}\ \right\}\right\}\\
\\
                     & \geqq &\inf\limits_{d\in L} \ \left\{\max\limits_{i \in X_2}\ \left\{\displaystyle \frac{-2\nabla p_i(x^*)^T d}{\gamma_i}\ \right\}\right\}\\
\\
     & \geqq & \inf\limits_{d\in L} \ \left\{\max\limits_{i \in X_2}\ \left\{\displaystyle \frac{-2\nabla p_i(x^*)^T d}{\gamma}\ \right\}\right\} =\inf\limits_{d\in L} \ \displaystyle \frac{F(d)}{-\gamma},
\end{array}$$
where $\displaystyle \gamma = \min\limits_{i\in I}\ \{\gamma_i\}$. Since $F(d)\geqq \rho$, we have $\beta \geqq \frac{\rho}{-\gamma}$. If 
$\Lambda^d =\emptyset$ for all $d\in \partial B(0,1)$, then $x^*$ is $\infty$-efficient, but if there exists $d$ such that 
$\Lambda^d \neq \emptyset$, then $x^*$ is $\beta$-efficient. As $\infty > \beta \geqq \frac{\rho}{-\gamma}$, the theorem holds. }}\\

Because $i\in X_2$, the passage from the first to the second inequality of the above demonstration is valid, and this also grants 
that $\rho>0$. Therefore $\frac{\rho}{-\gamma}>0$.

In order to limit the search space, another very useful result which gives a upper bound for the radius of efficiency of the 
solution $x^*\in$ \textit{Leff}$(VQFP')$ is presented.

{\tteo{Let $x^*\in$ \textit{Leff}$(VQFP')$ and $ M=\min\limits_{i\in X_1}\ \left\{\frac{ \left\|2\nabla p_i(x^*)\right\|}{\alpha}\right\}$. Suppose that 
$d^T \nabla^2 p_i(x^*) d \geqq \alpha > 0$, for some $d\in L$ and for all $i\in X_1 \neq \emptyset.$ If there does not exist another point $x'\in B(x^*,M)$ such that $\frac{f(x')}{g(x')} \leq \frac{f(x^*)}{g(x^*)}$, 
then $x^* \in$ \textit{Eff}$(VQFP')$.}\label{teo:BolaOndeNaoHaDominio}}\\

\begin{figure}[h]
\begin{center}
\begin{tikzpicture}[scale=0.7]

\draw [-triangle 45](0.5,-0.8) -- (0.5,8) node [left]{$x_2$} ;\draw [-triangle 45](0,-.5) -- (10,-.5) node [below]{$x_1$} ;

\draw [dashed, line width=.7pt](5,4) circle (4.15cm);
\draw [fill=gray!40, line width=.7pt](5,4) circle (3.5cm);
\draw [fill=white, line width=.7pt](5,4) circle (2cm);
\draw [dashed, line width=.7pt](5,4) circle (1.2cm);

\filldraw (5,4) node [left]{$x^*$} circle (.04cm);\filldraw (7.45,3) node [left]{$x'$} circle (.04cm);

\draw (5.2,5.4) node {\scriptsize $B(x^*,\frac{\rho}{-\gamma})$};\draw (6,6.15) node {\scriptsize $B(x^*,\beta)$};\draw (5.45,7.73) node {\scriptsize $B(x^*,P)$};\draw (6,8.35) node {\scriptsize $B(x^*,M)$};

\end{tikzpicture}

\caption{Some interesting neighborhoods of the solution $x^*\in$ \textit{Leff}$(VQFP')$}\label{fig:FronteirasDominancia}
\end{center}
\end{figure}

{\prova{We deduce from Lemma~\ref{lem:LambdaVazioNovaVersao} and Collorary~\ref{cor:LambdaNaoVazio} that if there exists a point $x'$ 
that dominates $x^*$ in direction $d$, then $X_2\neq \emptyset$ and $\lambda_{2}^d\leqq \left\|x'-x^*\right\| \leqq
\lambda_1^d$. If $P=\sup\limits_{d\in L} \ \{\lambda_1^d\}$, we have $\left\|x'-x^*\right\|\leqq P$. By the hypothesis, we also obtain
$$\begin{array}{rcl}
P &=& \sup\limits_{d\in L} \left\{ \min\limits_{i \in X_1} \left\{\displaystyle \frac{-2\nabla p_i(x^*)^T d}{d^T \nabla^2 p_i(x^*) d}\right\}\right\} \\
\\
                     & \leqq &\sup\limits_{d\in L} \left\{ \min\limits_{i \in X_1} \left\{\displaystyle \frac{-2\nabla p_i(x^*)^T d}{\alpha}\right\}\right\} \\
\\                     
                     & \leqq & \min\limits_{i \in X_1} \left\{\displaystyle \frac{\left\|2\nabla p_i(x^*)\right\|}{\alpha}\right\} = \ M.\\
                     
\end{array}$$
Therefore, if there does not exist a point that dominates $x^*$ in $B(x^*,M)$, then does not exist another poit in $\R^n$ with this property, and $x^* \in$
\textit{Eff}$(VQFP').$ }}\\

Because $i\in X_1$, the passage from the first to the second inequality of the above demonstration is valid. Theorem~\ref{teo:BolaOndeNaoHaDominio} only has usefulness if
$M<\infty.$

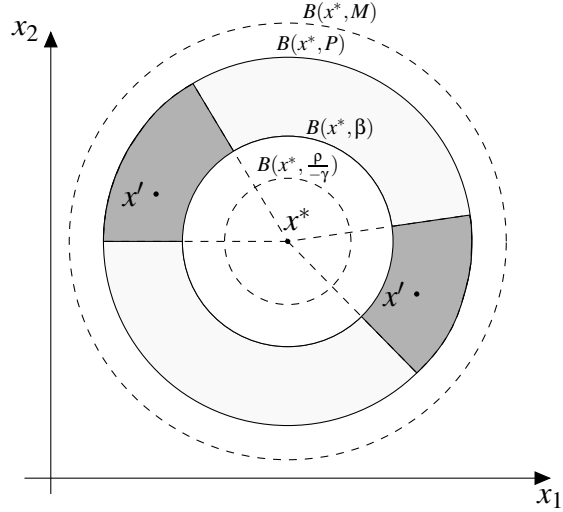
\begin{figure}[h]
\begin{center}
\begin{tikzpicture}[scale=0.7]

\draw [-triangle 45](0.5,-0.8) -- (0.5,8) node [left]{$x_2$} ;\draw [-triangle 45](0,-.5) -- (10,-.5) node [below]{$x_1$} ;

\draw [dashed](5,4) circle (4.15cm);
\draw [fill=gray!5](5,4) circle (3.5cm);
\draw [fill=white](5,4) circle (2cm);
\draw [dashed](5,4) circle (1.2cm);

\draw [fill=gray!60] (5,4) -- (7.45,1.5) .. controls (8.2,2.2) and (8.3,2.8) .. (8.43,3.3) -- (8.48,3.6) .. controls (8.5,4) and (8.5,4.4) .. (8.45,4.5) -- cycle;

\draw [fill=gray!60] (5,4) -- (1.5,4) .. controls (1.5,4.6) and (1.65,5.2) .. (2,5.8) -- (2,5.8) .. controls (2.27,6.2) and (2.5,6.6).. (3.2,7) -- cycle;

\draw [fill=white](5,4) circle (2cm);
\draw [dashed](5,4) circle (1.2cm);

\draw [dashed] (5,4) -- (7.45,1.5) .. controls (8.2,2.2) and (8.3,2.8) .. (8.43,3.3) -- (8.48,3.6) .. controls (8.5,4) and (8.5,4.4) .. (8.45,4.5) -- cycle;

\draw [dashed] (5,4) -- (1.5,4) .. controls (1.5,4.6) and (1.65,5.2) .. (2,5.8) -- (2,5.8) .. controls (2.27,6.2) and (2.5,6.6).. (3.2,7) -- cycle;

\filldraw (5,4) circle (.04cm);\draw (5.2,4) node [above]{$x^*$};\filldraw (7.45,3) node [left]{$x'$} circle (.04cm);\filldraw (2.5,4.9) node [left]{$x'$} circle (.04cm);

\draw (5.2,5.4) node {\scriptsize $B(x^*,\frac{\rho}{-\gamma})$};\draw (6,6.15) node {\scriptsize $B(x^*,\beta)$};\draw (5.45,7.73) node {\scriptsize $B(x^*,P)$};\draw (6,8.35) node {\scriptsize $B(x^*,M)$};

\end{tikzpicture}

\caption{Search space for $x'$ that dominates $x^*\in$ \textit{Leff}$(VQFP')$}\label{fig:FronteirasDominanciaReduzida}
\end{center}
\end{figure}

The problem of determining whether a locally Pareto-efficient optimal solution is dominated can be reduced by comparing it within a limited 
subset of $\mathbf{R}^n$ using the above results. Note that Collorary~\ref{cor:LambdaNaoVazio} provides the maximum radius $\beta=\inf\limits_{d\in L}\left\{\lambda_2^d\right\}$ of the neighborhood 
$B(x^*,\beta)$ where $x^*$ is not dominated. Collorary~\ref{cor:FMaiorRho} provides a lower bound for $\beta$, thus a new neighborhood 
$B(x^*,\frac{\rho}{-\gamma})$ for $x^*$. On the other hand, Theorem~\ref{teo:BolaOndeNaoHaDominio} provides two neighborhoods for $x^*$, 
one with radius $P=\sup\limits_{d\in L} \ \{\lambda_1^d\}$, and another with radius $M=\min\limits_{i\in X_1}\ \left\{\frac{\left\|2\nabla p_i(x^*)\right\|}{\alpha}\right\}$ 
as a upper bound for $P$. Therefore, we obtain four neighborhoods satisfying $B(x^*,\frac{\rho}{-\gamma})\subseteq B(x^*,\beta)
\subseteq B(x^*,P)\subseteq B(x^*,M)$.

The lower and upper bounds are more attractive computationaly to calculate. If we have a good computational search method to find a point 
that dominates the solution $x^*$, it is enought that this search is made in the subset $\bar{B}(x^*,P)\setminus B(x^*,\beta)$, or 
alternatively in the subset $\bar{B}(x^*,M)\setminus B(x^*,\frac{\rho}{-\gamma})$.

In general, the former neighborhoods are spheres in $\mathbf{R}^n$, but a particular case when the domain of the functions is 
$\mathbf{R}^2$ is illustrated in Figure~\ref{fig:FronteirasDominancia}. It shows four spheres centered in solution $x^*$: in dashed lines and 
nearer to the center appears the radius of neighborhood $B(x^*,\frac{\rho}{-\gamma})$; in continous lines and going little far away from 
the center appears the radius of neighborhoods $B(x^*,\beta)$ and $B(x^*,P)$, respectively; in shading appears the closed subset 
$\bar{B}(x^*,P)\setminus B(x^*,\beta)$; finally, in dashed lines appears the radius of the biggest neighborhood $B(x^*,M)$. If there exists a
point $x'$ that dominates $x^*$, it must belong to subset $\bar{B}(x^*,P)\setminus B(x^*,\beta)$, as shown in Figure~\ref{fig:FronteirasDominancia}.

Auxiliary computational methods can be designed to further reduce the size of subset $\bar{B}(x^*,P)\setminus B(x^*,\beta)$. 
It is enough to observe that $x'\in \bar{B}(x^*,P)\setminus B(x^*,\beta)$ if and only if there exists an unitary direction $d$ such that 
$\lambda_2^d\leqq \lambda \leqq \lambda_1^d$ and $x'=x^* + \lambda d$. Therefore, all those directions participating in subset 
$\bar{B}(x^*,P)\setminus B(x^*,\beta)$, but with $\Lambda^d=\emptyset$, can be excluded. Figure~\ref{fig:FronteirasDominanciaReduzida} 
illustrates this possibility, where there exists two shaded subsets in $\bar{B}(x^*,P)\setminus B(x^*,\beta)$, each containing a possible 
point $x'$. The shaded subsets represent only the directions $d$ which occurs $\lambda_2^d\leqq \lambda\leqq \lambda_1^d$.

\subsection{Radius of efficiency in the constrained case}\label{sec:comRestricoes}

Now, we can extend the results achieved in Section~\ref{sec:raioEmQFMP} for the constrained problem (VQFP) defined by 
feasible set $S$. Let $diam(S) = \sup \ \{\left\|x-y\right\|: x,y\in S\}$. The radius of efficiency $diam(S)$ is equivalent to 
the radius of efficiency $\infty$ in the (VQFP$'$). If $x^*$ is $\lambda$-efficient, then it is $\beta$-efficient, for all 
$\beta < \lambda\leqq diam(S)$. We say that $x^*$ is $diam(S)$-\textit{efficient} if it is globally Pareto-efficient optimal solution in $S$.

In (VOP), many constraint sets can be considered, however as a first contributed we consider this section the problem (VQFP) 
with the feasible set $S$ defined as $S=\{x\in \mathbf{R}^n\ |\ Cx\leqq b,\ C\in \R^{p\times n}, b\in \mathbf{R}^{p}\}$.

This choice makes it easy to calculate the distance from a point $x^*\in S$ to a boundary $S$, which is indispensable in the calculus of radius of efficiency. 
However, we believe that the results can extended for the more general constraint sets than linear inequalities.

To apply the concept radius of efficiency, we begin again from a solution $x^*\in$ \textit{Leff}$(VQFP)$. The objective 
is to calculate its radius of efficiency to validate if $x^*$ is globally Pareto-efficient optimal solution and also obtain some 
information regarding the neighborhoods of $x^*$. 

The next two results extend naturally the existing relation between the associated problems (VQFP) and (VQFP)$_{x^*}$ shown in 
Section~\ref{sec:raioEmQFMP}. Now we can admit the existence of a particular neighborhood of fixed radius for the solution 
$x^*$, where it is not dominated, and we can be calculate this radius using our metodology. The relationship is that the same 
fixed neighborhoood can be used in both problems.

{\tteo{Let $x^*\in$ \textit{Leff}$(VQFP)$. $x^*$ is $\lambda$-efficient for (VQFP) if and only if $x^*$ is $\lambda$-efficient for (VQFP)$_{x^*}$.}\label{teo:ProblemaAssociadoComRestricoes}}\\

{\prova{To calculate the radius of efficiency of the solution $x^*\in$ \textit{Leff}$(PMFQ)_{x^*}$,
we must find, if any exists, the values of $\lambda>0$ and a feasible unitary direction $d$ such that $x^*+\lambda d\in S$ dominates $x^*$. In other words, 
following the ideas until now, we must find a pair $(\bar{\lambda},\bar{d})$ which solves the following inequalities
\begin{eqnarray*}
f_i(x^*+\lambda d)-\frac{f_i(x^*)}{g_i(x^*)}g_i(x^*+\lambda d)\leqq f_i(x^*)-\frac{f_i(x^*)}{g_i(x^*)}g_i(x^*)=0,\ \mbox{for all}\ i\in I.
\end{eqnarray*}
As defined in the Section~\ref{sec:raioEmQFMP}, this implies the resolution of the inequalities~(\ref{eq:desig1}). Similarly, 
to calculate the radius of efficiency of the solution $x^*\in$ \textit{Leff}$(PMFQ)$,  we must find a pair $(\tilde{\lambda},\tilde{d})$ which 
solves the following inequalities
\begin{eqnarray*}
\frac{f_i(x^*+\lambda d)}{g_i(x^*+\lambda d)}\leqq \frac{f_i(x^*)}{g_i(x^*)} \ \ \Longleftrightarrow \ \ \lambda\left( \nabla p_i(x^*)^T d + \frac{\lambda}{2} d^T \nabla^2 p_i(x^*) d \right) \leqq 0,\ \mbox{for all}\ i\in I.
\end{eqnarray*}
This implies again the resolution of the inequalities~(\ref{eq:desig1}). Hence, if there exists a pair $(\bar{\lambda},\bar{d})$ for
the solution $x^*\in$ \textit{Leff}$(PMFQ)_{x^*}$, and a pair $(\tilde{\lambda},\tilde{d})$ for the same solution 
$x^*\in$ \textit{Leff}$(PMFQ)$, they are obtained from the same set of inequalities~(\ref{eq:desig1}). Therefore, 
$(\bar{\lambda},\bar{d})=(\tilde{\lambda},\tilde{d})$ and $x^*$ is $\lambda$-efficient for both problems, $\lambda=\bar{\lambda}=\tilde{\lambda}$. }}

{\tcor{$x^*\in$ \textit{Eff}$(PMFQ)$ if and only if $x^*\in$ \textit{Eff}$(PMFQ)_{x^*}$.}\label{cor:ProblemaAssociadoComRestricoes}}\\

{\prova{It is sufficient to replace the pair $(\lambda,\ d)$ to the pair $(\infty,\ d)$ in Theorem~\ref{teo:ProblemaAssociadoComRestricoes}.}}\\

We represent the active constraints set in $x^*\in S$ as $K(x^*)=\{j\in J\ |\ C_{j} x^*= b_j \}$, 
where $C_j$ is the $j$-th row of the matrix $C$. We also represent the tangent cone to $S$ in the point $x^*\in S$ as $\widetilde{T}(x^*)= \{y \in \R^n \ |\ C_{j} y \leqq 0,\ \forall j\in K(x^*),\ \|y\|=1\}$. 
We say that $T(x^*)=\widetilde{T}(x^*) \cap \partial B(0,1)$ is the \textit{feasible directions} set in the point $x^*$.

Next we extracted some results as in the unconstrained case, however by their similarities, some demonstrations are omitted. 

Given $x^*\in S$, for each $d\in T(x^*)$ such that $x^* + \lambda_\ell^d d \in S$, we say that the real number 
$\lambda_\ell^d \in\left(0,+\infty\right)$ is the \textit{limiting} for the radius of efficiency of $x^*$. And analogously to the 
Section~\ref{sec:raioEmQFMP} we define

$$\lambda_\ell^d = \min\limits_{j\in J,\ C_{j} d > 0} \left\{\frac{b_j - C_{j} x^*}{C_{j} d}\right\},$$

$$\lambda_{2}^d = \max\limits_{i \in X_2} \left\{\ \frac{-2\nabla p_i(x^*)^T d}{d^T \nabla^2 p_i(x^*) d}\ \right\},\ \ {\footnotesize \lambda_1^d = \left\{\begin{tabular}{ll} $ \min\left\{\min\limits_{i \in X_1} \left\{\displaystyle \frac{-2\nabla p_i(x^*)^T d}{d^T \nabla^2 p_i(x^*) d} \right\},\ \lambda_\ell^d \right\},$ & \mbox{if}\ \ $X_1\neq \emptyset$ \\
&\\
$\inf\left\{+\infty,\ \lambda_\ell^d \right\},$ & \mbox{if}\ \ $X_1 = \emptyset$, \end{tabular}\right.} $$

$$\Lambda_{2}^d=[\lambda_{2}^d,\ \infty),\qquad\qquad \Lambda_1^d= \left\{\begin{tabular}{ll}$(0,\ \lambda_1^d ],$ & \mbox{if}\ \ $X_1\neq \emptyset$ \\ $(0,\ \lambda_1^d ),$ & \mbox{if}\ \ $X_1=\emptyset$, \end{tabular}\right.\qquad \quad \Lambda^d = \Lambda_{2}^d \cap \Lambda_1^d,$$ 
where $p_i$, $i\in I$, $X_0$, $X_1$ and $X_2$ are the same functions and the same sets defined in Section~\ref{sec:raioEmQFMP}.

\begin{figure}[h]
\begin{center}
\begin{tikzpicture}

\draw [->](1,0.5) -- (1,5) node [above]{\small $r_i(\lambda)$};
\draw [->](.5,1) -- (7,1) node [right]{\small $\lambda$};

\draw [dashed](.9,2) -- (7,2) (0.5,2) node {{\tiny$p_{4} (x^*)$}}; \draw [line width=.7pt] (1,2) parabola[bend pos=0.5] bend +(0,4) +(4,-0.5);\draw [->, line width=1.2pt](4.9,2) -- (5.3,2);

\draw [dashed](.9,2.5) -- (7,2.5) (0.5,2.5) node {{\tiny$ \textcolor {black}{p_{3} (x^*)}$}};\draw[color=black, line width=.7pt] (1,2.5) sin (2.2,1.5) cos (4.5,5);\draw [color=black,<-, line width=1.2pt](3,2.5) -- (3.35,2.5);\draw [dashed](3.35,2.5) -- (3.35,.8) node[below] {{\scriptsize \textcolor{black}{$\lambda_1^{d}$}}};\draw [color=black,<-|, line width=1.2pt](2.6,1) -- (3.35,1);

\draw [dashed](.9,3) -- (7,3) (0.5,3) node {{\tiny$p_{2} (x^*)$}};\draw [line width=.7pt] (1,3) parabola[bend pos=0.5] bend +(0,2) +(5,-1.5);\draw [->, line width=1.2pt](5.2,3) -- (5.6,3);\draw [dashed](5.2,3) -- (5.2,0.8) node[below] {{\scriptsize $\lambda_{2}^{d}$}};\draw [|->, line width=1.2pt](5.2,1) -- (5.7,1);

\draw [dashed](.9,3.5) -- (7,3.5) (0.5,3.5) node {{\tiny$\textcolor {black}{p_{1} (x^*)}$}};\draw[color=black, line width=.7pt] (1,3.5) sin (2.5,2.5) cos (6,5.5);\draw [color=black,<-, line width=1.2pt](4,3.5) -- (4.4,3.5);

\draw (6.78,0.58) node {{\small $\Lambda^{d}=\emptyset $}};

\draw [blue, line width=.7pt](5.9,0.8) node [below]{{\scriptsize $\lambda_{\ell}^{d}$}} -- (5.9,5.7);

\draw [<-, color=black](4.65,5) -- (5.3,5.4);
\draw [->, color=black](5.3,5.4) -- (5.55,5);
\draw [color=black] (5.15,5.6) node {\small $1,3\in X_1$};

\draw [<-](2,4.56) -- (2.57,4.7);
\draw [->](2.57,4.7) -- (2.5,4.16);
\draw [] (3,4.9) node {\small $2,4\in X_2$};

\end{tikzpicture}
\caption{A search direction by $x'$ that dominates $x^*\in$ \textit{Leff}$(VQFP)_{x^*}$, $X_0=\emptyset$, $X_1=\{1,3\}$, $X_2=\{2,4\}$, 
$\lambda_1^{d}<\lambda_2^{d}\leqq\lambda_{\ell}^{d}$ and $\Lambda^{d}=\emptyset$}\label{fig:Direcao1BuscaComRestricoes}
\end{center}
\end{figure}

\begin{figure}[h]
\begin{center}
\begin{tikzpicture}

\draw [->](1,0.5) -- (1,5) node [above]{\small $r_i(\lambda)$};
\draw [->](.5,1) -- (7,1) node [right]{\small $\lambda$};

\draw [dashed](.9,2) -- (7,2) (0.5,2) node {{\tiny$p_{4} (x^*)$}}; \draw [line width=.7pt] (1,2) parabola[bend pos=0.5] bend +(0,4) +(4,-0.5);\draw [->, line width=1.2pt](4.9,2) -- (5.3,2);

\draw [dashed](.9,2.5) -- (7,2.5) (0.5,2.5) node {{\tiny$ p_{3} (x^*)$}};\draw [line width=.7pt] (1,2.5) parabola[bend pos=0.5] bend +(0,2) +(2,-1); \draw [->, line width=1.2pt](2.8,2.5) -- (3.2,2.5);

\draw [dashed](.9,3) -- (7,3) (0.5,3) node {{\tiny$p_{2} (x^*)$}};\draw [line width=.7pt] (1,3) parabola[bend pos=0.5] bend +(0,2) +(5,-1.5);\draw [->, line width=1.2pt](5.2,3) -- (5.6,3);\draw [dashed](5.2,3) -- (5.2,0.8) node[below] {{\scriptsize $\lambda_{2}^{d}$}};\draw [|->, line width=1.2pt](5.19,1) -- (6,1);\draw [color=blue,<-|, line width=1.2pt](3.8,1) -- (4.51,1);

\draw [dashed](.9,3.5) -- (7,3.5) (0.5,3.5) node {{\tiny$p_{1} (x^*)$}};\draw [line width=.7pt] (1,3.5) parabola[bend pos=0.5] bend +(0,2) +(3,-2);\draw [->, line width=1.2pt](3.35,3.5) -- (3.75,3.5);

\draw (6.45,0.6) node {{\small $\Lambda^{d}=\emptyset $}};\draw (6.45,4) node {{\small \textcolor{black}{$X_1=\emptyset $}}};

\draw [blue, line width=.7pt](4.5,0.8) node [below]{{\scriptsize $\lambda_{\ell}^{d}$}} -- (4.5,5.5) node [right]{{\scriptsize $\textcolor{black}{\lambda_1^{d}}\leftarrow \lambda_{\ell}^{d}$}} ;

\draw [] (2,1.2) node {\small$1,2,3,4\in X_2$};

\end{tikzpicture}
\caption{A search direction by $x'$ that dominates $x^*\in$ \textit{Leff}$(VQFP)_{x^*}$, $X_0=\emptyset$, $X_1=\emptyset$, $X_2=\{1,2,3,4\}$, 
$\lambda_{\ell}^{d}< \lambda_{2}^{d}$, $\lambda_1^{d}\leftarrow \lambda_{\ell}^{d}$ and $\Lambda^{d}=\emptyset$}\label{fig:Direcao4BuscaComRestricoes}
\end{center}
\end{figure}
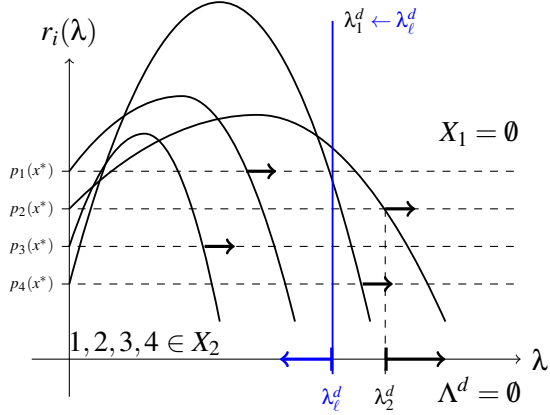

Let $x^*$ be locally Pareto-efficient optimal solution and $d\in T(x^*)$. Recall that if $X_0 = \emptyset$ and $X_1\neq \emptyset$, 
then $X_2=\emptyset$ cannot occurs, however, if $X_0 \neq \emptyset$ then $x^*$ is Pareto-efficient optimal solution in this direction. Also, if 
if $X_0 =\emptyset$, then either $X_1\neq \emptyset$ and $X_2\neq \emptyset$, or only $X_2\neq \emptyset$. Some of these possibilities are 
illustrated in Figures~\ref{fig:Direcao1BuscaComRestricoes},~\ref{fig:Direcao2BuscaComRestricoes},~\ref{fig:Direcao3BuscaComRestricoes}
and~\ref{fig:Direcao4BuscaComRestricoes}, including now $\lambda_{\ell}^d$ as possible limiting for the radius of efficiency in 
direction $d$. For each $i\in\{1,2,3,4\}$, the functions \begin{displaymath} r_i(\lambda) = p_i(x^*) + \lambda \nabla p_i(x^*)^T d + \frac{\lambda^2}{2} d^T \nabla^2 p_i(x^*) d, \end{displaymath} 
are plotted in the coordinates $(\lambda,r_i(\lambda))$, $\lambda>0$, whose graphics are parabolas and we can verify examples of 
the sets $X_0$, $X_1$, $X_2$ and $\Lambda^d$. Figure~\ref{fig:Direcao1BuscaComRestricoes} illustrates an example of direction where 
$X_1=\{1,3\}$ is formed by the indices of the functions $r_1$ and $r_3$, $X_2=\{2,4\}$ is formed by the indices of the fucntions $r_2$ and 
$r_4$, and $\Lambda^{d}=\emptyset$. In Figures~\ref{fig:Direcao1BuscaComRestricoes},~\ref{fig:Direcao2BuscaComRestricoes} and~\ref{fig:Direcao3BuscaComRestricoes} are
shown the cases in which $X_0=\emptyset$ and $X_1\neq \emptyset$, therefore $X_2\neq \emptyset$. However, 
Figure~\ref{fig:Direcao4BuscaComRestricoes} illustrates the case in which $X_1=\emptyset$ and only $X_2\neq \emptyset$.

\begin{figure}[h]
\begin{center}
\begin{tikzpicture}

\draw [->](1,0.5) -- (1,5) node [above]{\small $r_i(\lambda)$};
\draw [->](.5,1) -- (7,1) node [right]{\small $\lambda$};

\draw [dashed](.9,2) -- (7,2) (0.5,2) node {{\tiny$p_{4} (x^*)$}}; \draw [line width=.7pt] (1,2) parabola[bend pos=0.5] bend +(0,4) +(2,-0.8);\draw [->, line width=1.2pt](2.9,2) -- (3.3,2);\draw [dashed](2.9,2) -- (2.9,0.8) node[below] {{\scriptsize $\lambda_{2}^{d}$}};\draw [|->, line width=1.2pt](2.89,1) -- (3.5,1);

\draw [dashed](.9,2.5) -- (7,2.5) (0.5,2.5) node {{\tiny$ \textcolor {black}{p_{3} (x^*)}$}};\draw[color=black, line width=.7pt] (1,2.5) sin (3,1.5) cos (6,5);\draw [color=black,<-, line width=1.2pt](4.05,2.5) -- (4.45,2.5);\draw [dashed](4.45,2.5) -- (4.45,.8);\draw [color=blue,<-|, line width=1.2pt](1.66,1) -- (2.3,1);

\draw [dashed](.9,3) -- (7,3) (0.5,3) node {{\tiny$\textcolor {black}{p_{2} (x^*)}$}};\draw[color=black, line width=.7pt] (1,3) sin (3.5,1.3) cos (7,5.5);\draw [color=black,<-, line width=1.2pt](5.18,3) -- (5.58,3);

\draw [dashed](.9,3.5) -- (7,3.5) (0.5,3.5) node {{\tiny$p_{1} (x^*)$}};\draw [line width=.7pt] (1,3.5) parabola[bend pos=0.5] bend +(0,1.7) +(2.1,-2);\draw [->, line width=1.2pt](2.6,3.5) -- (3,3.5);

\draw (6.3,0.6) node {{\small $\Lambda^{d}=\emptyset $}};

\draw [blue, line width=.7pt](2.3,0.8) node[below] {{\scriptsize \textcolor{blue}{$\lambda_{\ell}^{d}$}}} -- (2.3,5.7) node [right]{{\scriptsize $\textcolor{black}{\lambda_1^{d}}\leftarrow \lambda_{\ell}^{d}$}};

\draw [<-, color=black](5.9,4.95) -- (5.8,5.3);
\draw [->, color=black](5.8,5.3) -- (6.68,5);
\draw [color=black] (6,5.5) node {\small $2,3\in X_1$};

\draw [<-](2.52,4.6) -- (2.8,4.6);
\draw [->](2.8,4.6) -- (2.4,4);
\draw [] (3.3,4.75) node {\small $1,4\in X_2$};

\end{tikzpicture}
\caption{A search direction by $x'$ that dominates $x^*\in$ \textit{Leff}$(VQFP)_{x^*}$, $X_0=\emptyset$, $X_1=\{2,3\}$, $X_2=\{1,4\}$, 
$\lambda_{\ell}^{d}< \lambda_{2}^{d}$, $\lambda_1^{d}\leftarrow \lambda_{\ell}^{d}$ and $\Lambda^{d}=\emptyset$}\label{fig:Direcao2BuscaComRestricoes}
\end{center}
\end{figure}

An important change is made in relation the parameter $\lambda_1^d$. Now, we have to include the possibility of 
$\lambda_1^d$ become $\lambda_{\ell}^d$. In the unconstrained problems we defined 
$\lambda_1^d=\min\limits_{i \in X_1} \left\{ \frac{-2\nabla p_i(x^*)^T d}{d^T \nabla^2 p_i(x^*) d}\right\}$ and, when 
$X_1\neq \emptyset$, the interest was in the directions $d$ such that $\lambda_2^d \leqq \lambda_1^d$. We continue interested in those 
directions, however observe that if occur $\lambda_{\ell}^d < \lambda_2^d$ does not make sense choose $\lambda$,
$\lambda_{\ell}^d < \lambda_2^d \leqq \lambda \leqq \lambda_1^d$ such that $x'=x^*+ \lambda d$ dominates $x^*$, because $x'\notin S$. 
Therefore, as was defined above, the correct is let the new $\lambda_1^d:=\min\{\lambda_1^d,\lambda_{\ell}^d\}$ when 
$X_1\neq \emptyset$, and $\lambda_1^d:=\inf\{+\infty,\lambda_{\ell}^d\}$ when $X_1=\emptyset$.

The role of the limiting $\lambda_{\ell}^d$ in the calculation of the radius of efficiency of $x^*\in$ \textit{Leff}$(VQFP)$ 
can be observed in the Figures~\ref{fig:Direcao1BuscaComRestricoes},~\ref{fig:Direcao2BuscaComRestricoes},~\ref{fig:Direcao3BuscaComRestricoes} 
and~\ref{fig:Direcao4BuscaComRestricoes}. $\lambda_{\ell}^d$ is represented by a continous vertical line that represents, in the position 
that it crosses the $\lambda$-axis of the cartesian plane and in the direction $d$, the boundary of the set $S$. 
Figure~\ref{fig:Direcao1BuscaComRestricoes} illustrates the case in which $\lambda_1^{d}<\lambda_2^{d}\leqq\lambda_{\ell}^{d}$, and 
therefore we have $\Lambda^{d}=\emptyset$. Figure~\ref{fig:Direcao3BuscaComRestricoes} illustrates the case in which 
$\lambda_2^{d}\leqq\lambda_1^{d}<\lambda_{\ell}^{d}$, and therefore we have $\Lambda^{d} \neq \emptyset$. On the other hand, 
Figure~\ref{fig:Direcao2BuscaComRestricoes} illustrates the case in which $\lambda_1^d$ become $\lambda_{\ell}^d$, that is 
$\lambda_1^{d}:=\lambda_{\ell}^{d}<\lambda_{2}^{d}$, and therefore $\Lambda^{d}=\emptyset$. Another possible case is illustrated in
Figure~\ref{fig:Direcao4BuscaComRestricoes}, in which the sets $X_0$ and $X_1$ are empty and we have only $X_2\neq \emptyset$. Therefore 
$\lambda_1^d$ become $\lambda_{\ell}^d$, that is $\lambda_1^{d}:=\lambda_{\ell}^{d}<\lambda_{2}^{d}$, and we obtain $\Lambda^{d}=\emptyset$.

\begin{figure}[h]
\begin{center}
\begin{tikzpicture}

\draw [->](1,0.5) -- (1,5) node [above]{\small $r_i(\lambda)$};
\draw [->](.5,1) -- (7,1) node [right]{\small $\lambda$};

\draw [dashed](.9,2) -- (7,2) (0.5,2) node {{\tiny$p_{4} (x^*)$}}; \draw [line width=.7pt] (1,2) parabola[bend pos=0.5] bend +(0,4) +(2,-0.8);\draw [->, line width=1.2pt](2.9,2) -- (3.3,2);\draw [dashed](2.9,2) -- (2.9,0.7) node[below] {{\scriptsize $\lambda_{2}^{d}$}};\draw [|->, line width=1.2pt](2.89,1) -- (3.5,1);

\draw [dashed](.9,2.5) -- (7,2.5) (0.5,2.5) node {{\tiny$ \textcolor {black}{p_{3} (x^*)}$}};\draw[color=black, line width=.7pt] (1,2.5) sin (3,1.5) cos (6,5);\draw [color=black,<-, line width=1.2pt](4.05,2.5) -- (4.45,2.5);\draw [dashed](4.45,2.5) -- (4.45,.7) node[below] {{\scriptsize \textcolor{black}{$\lambda_1^{d}$}}};\draw [color=black,<-|, line width=1.2pt](3.85,1) -- (4.46,1);

\draw [dashed](.9,3) -- (7,3) (0.5,3) node {{\tiny$\textcolor {black}{p_{2} (x^*)}$}};\draw[color=black, line width=.7pt] (1,3) sin (3.5,1.3) cos (7,5.5);\draw [color=black,<-, line width=1.2pt](5.18,3) -- (5.58,3);

\draw [dashed](.9,3.5) -- (7,3.5) (0.5,3.5) node {{\tiny$p_{1} (x^*)$}};\draw [line width=.7pt] (1,3.5) parabola[bend pos=0.5] bend +(0,1.7) +(2.1,-2);\draw [->, line width=1.2pt](2.6,3.5) -- (3,3.5);

\draw (6.5,0.5) node [color=black] {{\small $\Lambda^{d}\neq\emptyset $}};

\draw [blue, line width=.7pt](5.42,0.8) node [below]{{\scriptsize $\lambda_{\ell}^{d}$}} -- (5.42,5.5);
\draw [color=black,<->, line width=1.2pt] (2.89,.85) -- (4.46,.85);

\draw [<-, color=black](6.05,5.07) -- (6,5.3);
\draw [->, color=black](6,5.3) -- (6.67,5);
\draw [color=black] (6.23,5.45) node {\small $2,3\in X_1$};

\draw [<-](2.45,4.88) -- (2.8,5);
\draw [->](2.8,5) -- (2.3,4.1);
\draw [] (3.1,5.2) node {\small $1,4\in X_2$};

\end{tikzpicture}

\caption{A search direction by $x'$ that dominates $x^*\in$ \textit{Leff}$(VQFP)_{x^*}$, $X_0=\emptyset$, $X_1=\{2,3\}$, $X_2=\{1,4\}$, 
$\lambda_2^{d}\leqq\lambda_1^{d}<\lambda_{\ell}^{d}$ and $\Lambda^{d}\neq\emptyset$}\label{fig:Direcao3BuscaComRestricoes}
\end{center}
\end{figure}

Similar situations are possible, however among them we can exemplify two in which there exists a feasible point $x'$ that 
dominates $x^*$ in the direction $d\in T(x^*)$. First is shown in Figure~\ref{fig:Direcao3BuscaComRestricoes}, where there exists $\lambda>0$, 
$\lambda_2^{d}\leqq\lambda\leqq\lambda_1^{d} < \lambda_{\ell}^{d}$ such that $x^* + \lambda d$ dominates $x^*$, the interval 
$\Lambda^d$ is nonempty and $p(x')\leq p(x^*)$ has solution. Second can be observed in Figure~\ref{fig:Direcao4BuscaComRestricoes}, if 
we move $\lambda_{\ell}^d$ towards the right of the $\lambda$-axis until that $\lambda_2^{d}\leqq \lambda_{\ell}^{d}$, then $\lambda_1^d$ become 
$\lambda_{\ell}^d$, that is we would have $\lambda_2^{d}\leqq \lambda_{1}^{d}:=\lambda_{\ell}^d$ and obtain $\Lambda^d\neq \emptyset$.

{\tteo{Let $x^*\in$ \textit{Leff}$(VQFP)$. Then $x^*\in$ \textit{Eff}$(VQFP)$ if and only if for all $d\in T(x^*)$, $X_0\neq \emptyset$ or $\Lambda^d =\emptyset$.}\label{teo:LambdaVazioComRestricoes}}\\

{\prova{$(\Rightarrow)$ Identical to the shown in Theorem~\ref{teo:LambdaVazio}, considering now the directions in $T(x^*)$ and the fact 
of $S$ be limited or not in those directions.

$(\Leftarrow)$ By Theorem~\ref{teo:ProblemaAssociadoComRestricoes}, $x^*\in$\textit{Leff}$(VQFP)_{x^*}$. Let 
an arbitrary $d\in T(x^*)$ and suppose that $X_0 \neq \emptyset$ in direction $d$. Then there exists an index $i\in I$ such that one of 
items 1.(a) and 1.(b) or 2.(a) is satisfied, that is, $d^T \nabla^2 p_i(x^*) d>0$ and $\nabla p_i(x^*)^T d \geqq 0$, or 
$d^T \nabla^2 p_i(x^*) d\geqq 0$ and $\nabla p_i(x^*)^T d > 0$. Hence, $r_i(\lambda)$ grows indefinitely for $\lambda>0$ and given 
any neighborhood of $x^*$ in direction $d$, as large as it is, $p(x)\leq p(x^*)$ does not have solution for each $x=x^*+\lambda d$ in this 
neighborhood. Therefore, $x'=x^* + \lambda d$ does not dominates $x^*$ and $x^*$ is Pareto-efficient optimal solution in this direction. 
Figure~\ref{fig:Direcao3BuscaSemRestricoes} illustrates this posibility, where the function $r_1$, $1 \in X_0$, grows indefinitely for 
$\lambda>0$. On the other hand, suppose that $X_0 = \emptyset$ and $\Lambda^d=\emptyset$ in the direction $d$, then we can divide in two cases regarding the 
size of $X_1$. First, suppose also that $X_1=\emptyset$, then the Condition~\ref{cond:TudoNulo} grants that only $X_2\neq \emptyset$. Since 
$\Lambda^d=\emptyset$ imply $S$ is limited ($\lambda_{\ell}^d<+\infty$) in this direction, because $\lambda_1^d$ become $\lambda_{\ell}^d$, 
that is $\lambda_1^d:=\lambda_{\ell}^d < \lambda_2^d$ (see Figure~\ref{fig:Direcao4BuscaComRestricoes}), then for all
$\lambda \geqq \lambda_2^d$ any other point $x'=x^* + \lambda d\notin S$. Therefore, does not exist another feasible point that dominates $x^*$ 
and $x^*$ is Pareto-efficient optimal solution in this direction. On the second case, suppose that $X_1\neq\emptyset$ in the direction $d$, as $x^*$ is locally Pareto-efficient optimal solution, 
necesseraly $X_2\neq \emptyset$. Then we have $X_0= \emptyset$, $X_1\neq \emptyset$, $X_2\neq\emptyset$, $\Lambda^d=\emptyset$ and, 
as $\lambda_1^d:=\lambda_{\ell}^d$ whenever $\lambda_{\ell}^d\leqq \lambda_1^d$, we return to the unconstrained case as in 
Theorem~\ref{teo:LambdaVazio}. In fact, either $\lambda_1^d<\lambda_2^d$ and $p(x^*+\lambda d)\leq p(x^*)$ does not have solution for each $\lambda>0$, or by
Condition~\ref{cond:LambdaTudoIgual} and item 4.(a), $\lambda_2^d=\lambda_1^d=\frac{-2\nabla p_i(x^*)^T d}{d^T \nabla^2 p_i(x^*) d}$,
for all $i\in X_1\cup X_2$, and $p(x^*+\lambda d)\leq p(x^*)$ does not have solution for each $\lambda>0$. Therefore, does not exist another point 
$x'=x^*+\lambda d\in S$ that dominates $x^*$ and $x^*$ is Pareto-efficient optimal solution in this direction. Figures~\ref{fig:Direcao1BuscaComRestricoes} 
and~\ref{fig:Direcao2BuscaComRestricoes} illustrate this possibility, which can be divided in two cases. In the first case, 
Figure~\ref{fig:Direcao1BuscaComRestricoes} shows $\lambda_1^d=\min\limits_{i \in X_1} \left\{ \frac{-2\nabla p_i(x^*)^T d}{d^T \nabla^2 p_i(x^*) d} \right\} < \lambda_{\ell}^d$ 
and $\lambda_2^d \leqq \lambda_{\ell}^d$, however $\Lambda^d=\emptyset$ means that $\lambda_1^d < \lambda_2^d$ and 
inequalities~(\ref{eq:IdentidadeLambdaInicialMenor}), defined in Section~\ref{sec:raioEmQFMP}, are not satisfied for $\lambda>0$. Hence, does 
not exist another feasible point $x'=x^*+\lambda d$ that dominates $x^*$ in the direction $d$. In the second case, 
Figure~\ref{fig:Direcao2BuscaComRestricoes} shows $X_1\neq \emptyset$, $X_2\neq \emptyset$ and 
$\lambda_2^d \leqq \min\limits_{i \in X_1} \left\{\frac{-2\nabla p_i(x^*)^T d}{d^T \nabla^2 p_i(x^*) d} \right\}$, however 
$\Lambda^d=\emptyset$ means that there exists $j\in X_2$ such that $\lambda_\ell^d < \frac{-2\nabla p_j(x^*)^T d}{d^T \nabla^2 p_j(x^*) d}$ 
and $\lambda_\ell^d < \lambda_{2}^d$. Hence, for all $\lambda$ such that $\lambda_{2}^d \leqq \lambda \leqq \min\limits_{i \in X_1}\left\{ \frac{-2\nabla p_i(x^*)^T d}{d^T \nabla^2 p_i(x^*) d} \right\}$,
any other point $x'=x^*+\lambda d \notin S$, and does not exist another feasible point that dominates $x^*$ in this direction.

Therefore, given an arbitrary direction $d\in T(x^*)$, if $X_0\neq \emptyset$ we obtain $x^*$ non-dominated, or if 
$\Lambda^d=\emptyset$ we obtain again $x^*$ non-dominated. Hence, does not exist another point $x'=x^*+\lambda d\in S$ that dominates $x^*$, 
and $x^*\in$ \textit{Eff}$(VQFP)$. }}\\

It is important to observe that if $X_0 = \emptyset$, $X_1 = \emptyset$ and only $X_2\neq \emptyset$
in a direction $d\in T(x^*)$ it is possible to obtain $\Lambda^d=\emptyset$. However, in the unconstrained case always 
$\Lambda^{d}\neq\emptyset$ (see Figure~\ref{fig:Direcao4BuscaSemRestricoes}). Figure~\ref{fig:Direcao4BuscaComRestricoes} illustrates the 
constrained case, where we have $\lambda_{\ell}^d < \lambda_2^d$, in which $\lambda_1^d$ become $\lambda_{\ell}^d$, that is $\lambda_1^d:=\lambda_{\ell}^d < \lambda_2^d$, and 
we obtain $\Lambda^d=\emptyset$. As we define $\lambda_1^d=\inf\left\{+\infty,\ \lambda_\ell^d \right\}$, if for all $j\in X_2$, $\frac{-2\nabla p_j(x^*)^T d}{d^T \nabla^2 p_j(x^*) d}\leqq \lambda_\ell^d$, 
we return to the unconstrained case, because we would have $\lambda_2^d \leqq \lambda_1^d:=\lambda_{\ell}^d$ and $\Lambda^d \neq\emptyset$. That is, there exists
$\lambda \in \Lambda^d$ such that $x^*+\lambda d \in S$ dominates $x^*$ (see Figure~\ref{fig:Direcao4BuscaComRestricoes}, and move 
$\lambda_{\ell}^d$ towards the right of the $\lambda$-axis until $\lambda_2^{d}\leqq \lambda_{\ell}^{d}$).

The next results are extensions of the Colloraries~\ref{cor:LambdaNaoVazio},~\ref{cor:FMaiorRho} and 
Theorem~\ref{teo:BolaOndeNaoHaDominio}, so their proofs are equivalent and can be omitted. To simplify the presentation we define
$\hat{L}=\left\{d\in T(x^*)\ | \ \Lambda^d \neq \emptyset\right\}$.

{\tcor{Let $x^*\in$ \textit{Leff}$(VQFP)$ and $\beta = \inf\limits_{d\in \hat{L}}\ \left\{\lambda_2^d\right\}$. Then does not exist another 
point $x'\in B(x^*,\beta) \cap S$ such that $\frac{f(x')}{g(x')} \leq \frac{f(x^*)}{g(x^*)}$. }\label{cor:LambdaNaoVazioComRestricoes}}

{\tcor{Let $x^*\in$ \textit{Leff}$(VQFP)$ and $F(d) = \max\limits_{i\in X_2}\ \left\{ 2\nabla p_i(x^*)^T d \right\}.$ Suppose that exists 
$\rho\in \mathbf{R}$ such that for all $ d\in T(x^*)$ we have $F(d)\geqq \rho$. Then does not exist another point 
$x'\in B(x^*,\frac{\rho}{-\gamma}) \cap S$ such that $\frac{f(x')}{g(x')} \leq \frac{f(x^*)}{g(x^*)}$, where $\gamma < 0$, 
$\gamma = \min\limits_{i\in I} \left\{\gamma_i\right\}$ and $\gamma_i$ is the smallest negative eigenvalue of 
the matrix $\nabla^2 p_i(x^*)$, $i\in I$. }\label{cor:FMaiorRhoComRestricoes}}

{\tteo{Let $x^*\in$ \textit{Leff}$(VQFP)$ and $M=\min\left\{\min\limits_{i\in X_1} \frac{\left\|2\nabla p_i(x^*)\right\|}{\alpha}, \ diam(S)\right\}$. 
Suppose that $d^T \nabla^2 p_i(x^*) d \geqq \alpha > 0,$ for some $d\in \hat{L}$ and for all $i\in X_1\neq \emptyset$. If there does not 
exist another point $x'\in B(x^*,M)\cap S$ such that $\frac{f(x')}{g(x')} \leq \frac{f(x^*)}{g(x^*)}$, then 
$x^* \in$ \textit{Eff}$(VQFP)$. }\label{teo:BolaOndeNaoHaDominioComRestricoes}}\\

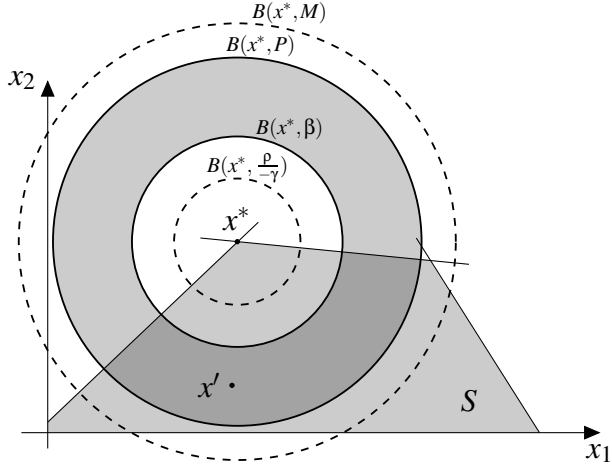
\begin{figure}[h]
\begin{center}
\begin{tikzpicture}[scale=0.7]


\draw [dashed, line width=.7pt](8.6,5.93) circle (4.15cm);
\draw [fill=gray!40, line width=.7pt](8.6,5.93) circle (3.5cm);
\draw [fill=white, line width=.7pt](8.6,5.93) circle (2cm);
\draw [dashed, line width=.7pt](8.6,5.93) circle (1.2cm);
\filldraw (8.6,5.93) node [above]{$x^*$} circle (.04cm);\filldraw (8.5,3.2) node [left]{$x'$} circle (.04cm);

\draw (8.8,7.33) node {\scriptsize $B(x^*,\frac{\rho}{-\gamma})$}; \draw (9.6,8.08) node {\scriptsize $B(x^*,\beta)$}; \draw (9.05,9.66) node {\scriptsize $B(x^*,P)$}; \draw (9.6,10.28) node {\scriptsize $B(x^*,M)$};

\draw [-triangle 45](5,2) -- (5,9) node [left]{$x_2$}; \draw [-triangle 45](4.5,2.3) -- (15.5,2.3) node [below]{$x_1$} ;
\draw (5,2.5)--(9,6.3);\draw (7.9,6)--(13,5.5);\draw (12,6)--(14.34,2.3);
\draw (13,2.9) node {$S$};
\fill[opacity=0.2] (5,2.3)--(5,2.5)--(8.6,5.93)--(12.25,5.58)--(14.34,2.3);

\end{tikzpicture}

\caption{Some interesting neighborhoods of the solution $x^*\in$ \textit{Leff}$(VQFP)$}\label{fig:FronteirasDominanciaComRetricoes}
\end{center}
\end{figure}

Similarly to the shown in Section~\ref{sec:raioEmQFMP}, we can identify in the Colloraries~\ref{cor:LambdaNaoVazioComRestricoes},~\ref{cor:FMaiorRhoComRestricoes} and 
Theorem~\ref{teo:BolaOndeNaoHaDominioComRestricoes} four important subsets related to solution $x^*\in$ \textit{Leff}$(PMFQ)$ satisfying 
$(B(x^*,\frac{\rho}{-\gamma})\cap S) \subseteq \left(B(x^*,\beta) \cap S \right) \subseteq \left(B(x^*,P) \cap S \right) \subseteq \left(B(x^*,M) \cap S\right)$, 
where $P=\sup\limits_{d\in \hat{L}} \ \{\lambda_1^d\}$. If we have a good computational search method to find a point $x'$ that dominates 
the solution $x^*$, it is enough that this search is made in the subset $\left(\bar{B}(x^*,P)\cap S\right)\setminus \left(B(x^*,\beta)\cap S\right)$, or 
alternatively in the subset $\left(\bar{B}(x^*,M)\cap S\right) \setminus (B(x^*,\frac{\rho}{-\gamma})\cap S)$. This way, we can present the 
following collorary.

{\tcor{Let $x^*\in$ \textit{Leff}$(PMFQ)$. Suppose that the hypothesis of 
Corollaries~\ref{cor:LambdaNaoVazioComRestricoes},~\ref{cor:FMaiorRhoComRestricoes} and 
Theorem~\ref{teo:BolaOndeNaoHaDominioComRestricoes} are satisfied, if there exists $x'$ that dominates $x^*$ in problem (VQFP), 
then $x'\in \left(\bar{B}(x^*,P)\cap S\right)\setminus \left(B(x^*,\beta)\cap S\right)$. }\label{cor:EnglobaTudo}}\\

Figure~\ref{fig:FronteirasDominanciaComRetricoes} explains the Corollary~\ref{cor:EnglobaTudo}, it is shown a
solution $x^*\in$ \textit{Leff}$(PMFQ)$ in the boundary of the $S\subseteq \mathbf{R}^2$. In dashed lines, two subsets of interest are 
shown: $B(x^*,\frac{\rho}{-\gamma})\cap S$ and $B(x^*,M)\cap S$. And in continous lines, another two are shown: $B(x^*,\beta)\cap S$ and $B(x^*,P)\cap S$. 
If there exists a point $x'$ that dominates $x^*$, it must belong to subset $\left(\bar{B}(x^*,P)\cap S\right)\setminus \left(B(x^*,\beta)\cap S\right)$.

\section{Conclusions}\label{cap:Conclusions}

The main contribution of this work is the development of necessary and sufficient Pareto optimality conditions for the solutions of a 
particular vector optimization problem, where each objective function consists of a ratio of two quadratic functions and the feasible 
set is defined by linear inequalities. We introduce the new concept of radius of efficiency in order to identify the neighborhoods of a 
locally Pareto-efficient optimal solution in vector optimization problems. We show how to calculate the radius of the very useful two 
spherical regions centered in this solution, and if there exists another point that dominates the former solution, it belongs to the 
subtraction of those spherical regions. In this process we may conclude that the solution is also globally optimal. Theorems are 
established for it. These results might be useful to determine termination criteria in the development of algorithms, and new extensions 
can be established from these to more general vector optimization problems in which quadratic approximations are used locally. In future 
work we plan to develop algorithms using the concept presented here.

\section*{Acknowledgements}
This work was partially supported by Coordination for the Improvement of Higher Level Personnel of Brazil (CAPES), Fund to Support Teaching, 
Research and Extension of Unicamp (FAEPEX) under grant 534/09 and Spain’s Ministry of Science and Technology under grant MTM2007-63432. 

\bibliographystyle{plain}
\bibliography{bibliografiaTese}

\end{document}